\documentclass[a4paper]{article}
\def\sup{{\rm sup}}
\def\max{{\rm max}}
\def\min{{\rm min}}
\def\htt{{\rm ht}}
\def\Max{{\rm Max}}
\def\Spec{{\rm Spec}}
\def\Min{{\rm Min}}
\begin{document}
\baselineskip=14pt
\newtheorem{thm}{Theorem}[section]
\newtheorem{example}[thm]{Example}
\newtheorem{remark}[thm]{Remark}
\newtheorem{pro}[thm]{Proposition}
\newtheorem{lem}[thm]{Lemma}
\newtheorem{cor}[thm]{Corollary}
\renewcommand{\thefootnote}{\fnsymbol{footnote}} 
\renewcommand{\thesection}{\arabic{section}.}
\renewcommand{\thethm}{\thesection\arabic{thm}}
\begin{center}
{\small Journal of Pure and Applied Algebra \tt (To Appear)} \bigskip

{\Large\bf On the prime ideal structure of tensor products of
  algebras}
\vspace{1cm}

{\large Samir Bouchiba$^a$, David E. Dobbs$^b$, Salah-Eddine Kabbaj$^{c,}$\footnote{Research partially supported by the Arab Fund for Economic and Social Development.}}
\bigskip\bigskip

{\small $^a$ Department of Mathematics, University Moulay Ismail,
Meknes 50000, Morocco \medskip

$^b$ Department of Mathematics, University of Tennessee,
Knoxville, TN 37996, USA \smallskip

$^c$ Department of Mathematics, Harvard University,
Cambridge, MA 02138, USA}
\end{center}
\bigskip

\hrule\medskip
\noindent{\large\bf Abstract} \medskip

\noindent This paper is concerned with the prime spectrum of a tensor 
product of
algebras over a field. It seeks necessary and
sufficient conditions
for such a tensor product to have the S-property, strong S-
property, and
catenarity. Its main results lead to new examples of stably 
strong S-rings
and universally catenarian rings. The work begins by 
investigating the
minimal
prime ideal structure.
Throughout, several results on polynomial rings are 
recovered, and
numerous examples are provided to illustrate the scope and 
sharpness of the
results.
\medskip

\medskip


\medskip\hrule\bigskip

\begin{section}{Introduction}

All rings and algebras considered in this paper are 
commutative
with identity element and, unless otherwise specified, are 
assumed to be
non-zero. All ring homomorphisms are unital. Throughout, 
$k$ denotes a
field. We shall
use $t.d.(A:k)$, or $t.d.(A)$ when no confusion is likely, 
to denote the
transcendence degree of a $k$-algebra $A$ over $k$ (for 
nondomains,
$t.d.(A)=\sup\{t.d.(\frac Ap):p\in \Spec(A)\}$), and $k_A(p)$ 
to denote
the
quotient field of $\frac Ap$, for each prime ideal $p$ of 
$A$. Also, we use
$\Spec(A)$, $\Max(A)$, and $\Min(A)$ to denote the sets of 
prime ideals,
maximal ideals, and minimal prime ideals, respectively, of 
a ring $A$, and
$\subset$
to denote proper inclusion. Recall that an integral domain 
$A$ of finite
Krull dimension $n$ is a Jaffard domain if its valuative 
dimension,
$\dim_v(A)$,
is also $n$. A locally Jaffard domain is a finite-
dimensional domain $A$
such that $A_p$ is a Jaffard domain for each $p\in \Spec(A)$.
Finite-dimensional Pr\" {u}fer domains and Noetherian 
domains are locally
Jaffard domains. We
assume familiarity with the above concepts, as in [1] and 
[15]. Any
unreferenced
material is standard, as in [12], [18], and [20].

Since the EGA of Grothendieck [13], a few works in the 
literature have
explored
the prime ideal structure of tensor products of
$k$-algebras (cf. [23], [24], [26], [3], and [4]). These 
have
mainly been
concerned with dimension theory in specific contexts, such 
as tensor
products
of fields,
AF-domains, or pullbacks. At present, the general situation 
remains
unresolved. By analogy with known studies on
polynomial rings, the investigation of some chain
conditions may be expected to cast light on the spectrum of 
such
constructions.
Thus, we focus here on an in-depth study of central notions 
such as the
S-property, strong
S-property, and
catenarity. In particular, our main result, Theorem 4.13, 
allows us to
provide new families of stably strong S-rings and
universally catenarian rings. Throughout, several results on
polynomial rings are recovered and
numerous examples are provided to illustrate the scope and 
sharpness of the
main results.

In order to treat Noetherian domains and
Pr\"{u}fer domains in a unified manner, Kaplansky [18] 
introduced the
concepts of S(eidenberg)-domain and strong S-ring. A domain 
$A$ is called
an S-domain if, for each height-one
prime ideal $p$ of $A$, the extension $pA[X]$ to the 
polynomial
ring in one variable also has height $1$. A commutative 
ring $A$ is said
to be a strong S-ring if $\frac Ap$ is an S-domain for each 
$p\in \Spec(A)$.
It is noteworthy that while $A[X]$ is always an S-domain 
for any domain
$A$ [11], $A[X]$ need not be a strong S-ring even when $A$
is a strong S-ring. Thus, as in [19], $A$ is said to be a 
stably
strong S-ring (also called a universally strong S-ring ) if 
the polynomial
ring $A[X_1,...,X_n]$ is a strong S-ring for each
positive integer $n$. The study of this class of rings was 
initiated by
Malik and Mott
[19] and further developed in [16] and [17]. An example of 
a strong
S-domain which
is not a stably strong S-domain was constructed in [8].

As in [5], we say that a domain $A$ is catenarian if $A$ is 
locally
finite-dimensional (LFD for short) and, for each pair 
$P\subset Q$ of
adjacent prime
ideals of $A$, $\htt(Q)=1+\htt(P)$; equivalently, if for
any prime ideals $P\subseteq Q$ of $A$, all the saturated
chains in $\Spec(A)$ between $P$ and $Q$ have the same 
finite length. Note
that catenarity is not stable under adjunction of 
indeterminates.
Thus, as in [5], a domain $A$ is said to be universally 
catenarian if
$A[X_1,...,X_n]$
is catenarian for each positive integer $n$. Cohen-Macaulay 
domains [20] or
LFD Pr\"ufer domains [7] are universally
catenarian; and so are
domains of valuative dimension 1 [5] and LFD domains of 
global dimension 2
[6].
Finally, recall that any universally catenarian domain is a 
stably strong
S-domain
[5, Theorem 2.4].

In Section 2, we extend the definitions of the S-property 
and catenarity to
the context of
arbitrary rings
(i.e., not necessarily domains).
Section 3 investigates the minimal prime
ideal structure in tensor products of $k$-algebras. Vamos 
[25] proved that
if $K$
and $L$ are field extensions of $k$, then the minimal 
prime ideals
of $K\otimes_kL$ are pairwise comaximal. We give an example 
to show that
this result fails for arbitrary domains $A$ and $B$ that are
$k$-algebras, and then show
that the minimal prime ideals of $A\otimes_kB$ are pairwise 
comaximal
provided that
$A$ and $B$ are integrally closed domains. As an 
application, we establish
necessary and
sufficient conditions for $A\otimes_kB$ to be an S-ring, 
and therefore
extend (in Theorem 3.9) the
known result that $A[X_1,...,X_n]$
is an S-domain for any domain $A$ and any integer $n\geq 1$ 
[11,
Proposition 2.1].
Our purpose in Section 4 is to study conditions under which 
tensor
product
preserves the strong S-property and catenarity. We begin 
with a result of
independent interest (Proposition 4.1) characterizing the LFD 
property for
$A\otimes_kB$. Also noteworthy is Corollary 4.10 stating 
that the tensor
product of two field extensions of $k$, at least one of 
which is of finite
transcendence degree, is universally catenarian.
Our main theorem (4.13) asserts that: given an LFD $k$-
algebra $A$ and an
extension
field $K$ of $k$ such that either $t.d.(A:k)< \infty $ or 
$t.d.(K:k)< \infty
$, let $B$ be a transcendence basis of $K$ over $k$ and $L$ 
be the separable
algebraic closure of $k(B)$ in $K$, and assume that $[L:k
(B)]<
\infty$; then if $A$ is a stably strong S-ring (resp., 
universally
catenarian and the minimal prime ideals of $K\otimes _kA$ 
are pairwise
comaximal),
$K\otimes _kA$ is a stably
strong S-ring (resp., universally catenarian). This result 
leads to new
families of stably
strong S-rings and universally catenarian rings. Section 5 
displays examples
illustrating
the limits of the results of earlier sections. The section 
closes with an
example of a discrete rank-one
valuation domain $V$ (hence universally catenarian) such 
that $V\otimes_kV$
is not
catenarian, illustrating the importance of assuming $K$ is 
a field in
Theorem 4.13.
\end{section}\bigskip

\begin{section} {Preliminaries}

In this section, we extend the notions of S-domain and 
catenarian domain
to the context of arbitrary rings (i.e., not necessarily 
domains). We then
state some
elementary results and recall certain basic facts about 
tensor
products of $k$-algebras, providing a suitable background 
to the rest of the
paper.\bigskip

Consider the following four properties that a ring $A$ may 
satisfy:\\
$(P_1): \frac AP$ is an S-domain for each $P\in \Min(A)$.\\
$(P_2): \htt(P)=1 \Rightarrow \htt(P[X])=1$, for each
$P\in \Spec(A)$.\\
$(Q_1): A$ is LFD and $\htt(Q)=1+\htt(P)$ for each pair 
$P\subset Q$ of adjacent
prime

\hspace{6.85mm}ideals of $A$.\\
$(Q_2): \frac AP$ is a catenarian domain for each $P\in \Min
(A)$.\bigskip

It is clear that a domain $A$ satisfies $(P_1)$ (resp.,
$(P_2)$) if and only
if $A$ is an S-domain; and that a domain $A$ satisfies 
$(Q_1)$ (resp.,
$(Q_2)$) if and only if $A$ is catenarian. Some of these 
observations carry
over to arbitrary rings. Using the basic facts from [18, p. 
25], we verify
easily that $(P_1)\Rightarrow (P_2)$; and that $(Q_1)
\Rightarrow
(Q_2)$. However,
the inverse implications do not hold in general. The next 
example
illustrates
this fact.
\begin{example} There exists a (locally) 
finite-dimensional
ring
$A$ which satisfies both $(P_2)$ and $(Q_2)$ but neither 
$(P_1)$
nor $(Q_1)$.\end{example}

Let $V:=k( X)[Y]_{( Y)}=k(
X)+m$, where $m:=YV$. Let $R=k+m$. There exist two 
saturated chains
in $\Spec(R[Z])$ of the form:
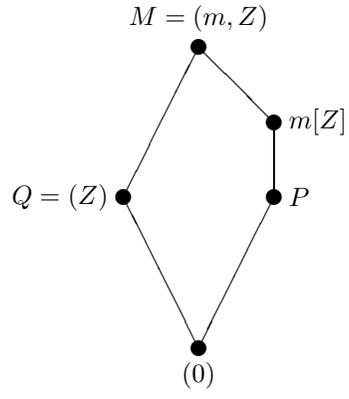
\[\setlength{\unitlength}{1mm}
\begin{picture}(40,56)(10,-8)
\put(20,0){\line(-1,2){10}}
\put(10,20){\line(1,2){10}}
\put(20,40){\line(1,-1){10}}
\put(30,30){\line(0,-1){10}}
\put(30,20){\line(-1,-2){10}}
\put(20,40){\circle*{2}}
\put(20,42){\makebox(0,0)[b]{$M=(m,Z)$}}
\put(10,20){\circle*{2}}
\put(8,20){\makebox(0,0)[r]{$Q=(Z)$}}
\put(20,0){\circle*{2}}
\put(20,-2){\makebox(0,0)[t]{$(0)$}}
\put(30,20){\circle*{2}}
\put(32,20){\makebox(0,0)[l]{$P$}}
\put(30,30){\circle*{2}}
\put(32,30){\makebox(0,0)[l]{$m[Z]$}}
\end{picture}\]

\noindent Indeed, let $I=PQ$ and $A=(\frac{R[Z]}I)_{\frac 
MI}$. Then $A$ is
a two-dimensional
quasilocal ring, and hence trivially satisfies $(Q_2)$. 
Further, part of
$\Spec(A)$ displays as follows:
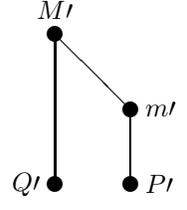
\[\setlength{\unitlength}{1mm}
\begin{picture}(30,30)(0,-5)
\put(0,0){\line(0,2){20}}
\put(0,20){\line(1,-1){10}}
\put(10,10){\line(0,-1){10}}
\put(0,0){\circle*{2}}
\put(-2,0){\makebox(0,0)[r]{$Q\prime$}}
\put(0,20){\circle*{2}}
\put(0,22){\makebox(0,0)[b]{$M\prime$}}
\put(10,10){\circle*{2}}
\put(12,10){\makebox(0,0)[l]{$m\prime$}}
\put(10,0){\circle*{2}}
\put(12,0){\makebox(0,0)[l]{$P\prime$}}
\end{picture}\]
\noindent where $M^{\prime }=(\frac MI)_{\frac MI}$, $Q^
{\prime }=
(\frac QI)_{\frac MI}$, $m^{\prime }=(\frac{m[Z]}I)_{\frac 
MI}$, and $
P^{\prime }=(\frac PI)_{\frac MI}$. It is clear that $m^
{\prime}$ is the
unique prime ideal of $A$ of height $1$. By [8, Example 5],
$\htt(m[Z])=\htt(m[Z,T])=2$, so that $\htt(m^{\prime}[T])=1$.
Thus $A$ satisfies $(P_2)$. Now, $\frac
A{Q^{\prime }}\cong (\frac{R[Z]}Q)_{\frac MQ}\cong R$
is not an S-domain, since $\htt(m)=1$ and $\htt(m[Z])=2$, 
whence $A$ does not
satisfy $(P_1)$. Moreover, $A$ fails to satisfy $(Q_1)$,
since $Q^{\prime}\subset M^{\prime}$ is a saturated chain 
in $\Spec(A)$
such that $\htt(M^{\prime })=2\neq
1+\htt(Q^{\prime })=1$. $\diamondsuit$\bigskip

By avoiding a feature of Example 2.1, we shall find a 
natural context in
which
$(P_2)$ implies $(P_1)$, and $(Q_2)$ implies $(Q_1)$. Let 
us say
that a ring $A$ satisfies MPC (for Minimal Primes 
Comaximality) if the
minimal prime ideals in $A$ are pairwise comaximal; i.e., 
if each maximal
ideal of $A$ contains only one minimal prime ideal. In the 
literature, MPC
has
also been termed ``locally irreducible", presumably because 
any domain
evidently satisfies MPC.
\begin{remark} Let $A$ be a ring satisfying 
MPC. Then:

\noindent a) $A$ satisfies $(P_1)$ (resp., $Q_1$) if and 
only if $A$
satisfies $(P_2)$ (resp., $Q_2$).

\noindent b) $S^{-1}A$ satisfies MPC for any multiplicative 
subset $S$ of
$A$.

\noindent c) $A[X_1,...,X_n]$ satisfies MPC for all 
integers $n\geq 1$.\end{remark}

\noindent {\bf Proof.} The proof of (a) may be left to the 
reader.
Now, (b) follows from basic facts about localization, while 
(c) is
immediate since the minimal prime ideals of $A[X_1,...,X_n]
$ are of
the form $p[X_1,...,X_n]$, where $p\in \Min(A)$.\bigskip

We now extend the domain-theoretic definitions of the S-
property and
catenarity
to the MPC context. A ring $A$ is called an S-ring if it 
satisfies
MPC and $(P_1)$; equivalently, MPC and $(P_2)$. A ring $A$ 
is said to be
catenarian
if $A$ satisfies MPC and $(Q_1)$; equivalently,
MPC and $(Q_2)$. It is useful to note that if $A$ is an S-
ring
(resp., a catenarian ring), then so is $A_S$ ($=S^{-1}A$), 
for each
multiplicative subset $S$ of $A$.\bigskip

Next, we extend a domain-theoretic result of Malik and Mott 
[19, Theorem
4.6].
\begin{pro} Let $A\subseteq T$ be an 
integral ring
extension. If $T$ is a strong S-ring (resp., stably strong 
S-ring), then
so is $A$.\end{pro}

\noindent {\bf Proof.} Let $p\in \Spec(A)$. Since $T$ is an 
integral
extension of $A$, the Lying-over Theorem provides $P\in Spec
(T)$ such
that $P\cap A=p$. Hence $\frac TP$ is an integral extension 
of $\frac Ap$,
and $\frac TP$ is a strong S-domain by hypothesis. 
Consequently, by
[19, Theorem 4.6], $\frac Ap$ is a (strong) S-domain. The 
``stably
strong S-ring" assertion follows from the ``strong S-ring" 
assertion
since $A[X_1,...,X_n]\subseteq T[X_1,...,X_n]$ inherits 
integrality
from $A\subseteq T$. $\diamondsuit$\bigskip

Proposition 2.5 generalizes the following domain-theoretic 
result.
\begin{pro} [5, Corollary 6.3] Let $A$ 
be a
one-dimensional
domain. Then the following conditions are equivalent:

i) $A$ is universally catenarian;

ii) $A[X]$ is catenarian;

iii) $A$ is a stably strong S-domain;

iv) $A$ is a strong S-domain;

v) $A$ is an S-domain.$\diamondsuit$\end{pro}
\begin{pro} Let $A$ be a one-
dimensional ring. Then\\
a) The following three conditions are equivalent:

i) $A$ is a stably strong S-ring;

ii) $A$ is a strong S-ring;

iii) $A$ satisfies $(P_2)$.

\noindent b) Suppose, in addition, that $A$ satisfies MPC. 
Then
(i)-(iii) are equivalent to each of (iv)-(vi):

iv) $A$ is universally catenarian;

v) $A[X]$ is catenarian;

vi) $A$ is an S-ring.\end{pro}

\noindent {\bf Proof.} a) It is trivial that
$(i)\Rightarrow (ii)$(even without one-dimensionality). 
Also, any
field is an S-domain. As $\dim(A)=1$, (ii) is therefore 
equivalent to
the requirement that $\frac AQ$ is an S-domain for each 
$Q\in \Min(A)$. This
requirement is obviously equivalent to (iii).
Thus, $(ii)\Leftrightarrow (iii)$.\\
\noindent $ii)\Rightarrow i)$ Clearly, it suffices to prove 
that $\frac Ap$
is a stably strong S-domain for each $p\in \Min(A)$. By
Proposition 2.4, this assertion holds, since for any $p\in 
\Min(A)$, $\frac
Ap$ is either a field or a one-dimensional strong
S-domain.\\
b) $(iii)\Leftrightarrow (P_1)\Leftrightarrow (vi)$, since 
$A$ satisfies
MPC.\\
\noindent $v)\Rightarrow vi)$ Let $p\in \Min(A)$. Then
$\frac Ap[X]\cong \frac{A[X]}{p[X]}$ is a catenarian 
domain, since $p[X]\in
\Min(A[X])$. So, by
Proposition 2.4, $\frac Ap$ is an S-domain. Hence (in view 
of the MPC
condition), $A$ is an S-ring.\\
\noindent $vi)\Rightarrow iv)$ It suffices to prove that 
$\frac Ap$ is
universally
catenarian, for each minimal prime ideal $p$ of $A$. This 
holds by
Proposition 2.4, since for any $p\in \Min(A)$, $\frac Ap$
is either a field or a one-dimensional S-domain. The proof 
is complete.
$\diamondsuit$\bigskip

For the convenience of the reader, we close this section by 
discussing some
basic facts connected with the tensor product of $k$-
algebras. These will
be used frequently in the sequel without explicit mention.\bigskip

Let $A$ and $B$ be two $k$-algebras. If $A^{\prime }$ is an 
integral
extension of $A$, then $A^{\prime}\otimes _kB$ is an 
integral extension of
$A\otimes _kB$.
If $S_1$ and $S_2$ are multiplicative subsets of $A$ and 
$B$,
respectively, then $S_1^{-1}A\otimes _kS_2^{-1}B\cong S^{-1}
(A\otimes _kB)$,
where $S:=\{s_1\otimes s_2:s_1\in S_1$ and $s_2\in S_2\}$. 
Recall also that
if $A$ is an integral domain, then $\htt(p)+t.d.(\frac Ap)
\leq t.d.(A)$, for
each $p\in \Spec(A)$
(cf. [21, p. 37] and [28, p. 10]). It follows that $\dim(A)\leq 
t.d.(A)$ for any ring $A$.
Moreover, we assume familiarity with the natural 
isomorphisms for tensor
products.
In particular, we identify $A$ and $B$ with their canonical
images in $A\otimes_kB$. Also, $A\otimes_kB$ is a free 
(hence faithfully
flat)
extension of $A$ and $B$. Here we recall that if $R\hookrightarrow S$ is a flat ring extension and $P\in \Min(S)$, then $P\cap R \in \Min(R)$ by going-down. Finally, we refer the reader to 
the useful result
of
Wadsworth [26, Proposition 2.3] which yields a 
classification of the prime
ideals of $A\otimes_kB$ according to their contractions to 
$A$ and $B$.
\end{section}\bigskip

\begin{section} {Minimal prime ideal structure\\ and S-
property}

This section studies the transfer of the MPC property and S-
property to
tensor products of $k$-algebras. As a prelude to this, we 
first
investigate the minimal prime ideal structure of such 
constructions.
In [25, Corollary 4], Vamos proved that if $K$ and $L$ are 
field extensions
of $k$, then
$K\otimes_kL$ satisfies MPC. We first illustrate
by an example the failure of this result for arbitrary
$k$-algebras $A$ and $B$, and then show that $A\otimes_kB$ 
satisfies MPC
provided
$A$ and $B$ are integrally closed domains. As an 
application, we establish
necessary and
sufficient conditions for $A\otimes_kB$ to be an S-ring, 
and therefore
extend the
known result that $A[X_1,...,X_n]$
is an S-domain, for any domain $A$ and any integer
$n\geq 1$ [11, Proposition 2.1]. Throughout Sections 3 and 
4, LO (resp., GD)
refers to the
condition ``Lying-over'' (resp., ``Going-down''), as in [18, 
p. 28].\bigskip

We begin by providing a necessary condition for 
$A\otimes_kB$ to satisfy
MPC.
\begin{pro} a) If $C\subseteq D$ is a 
ring extension
satisfying LO and GD, and $D$ satisfies MPC, then $C$ 
satisfies MPC.

\noindent b) If $A$ and $B$ are $k$-algebras such that 
$A\otimes _kB$
satisfies MPC,
then $A$ and $B$ each satisfy MPC.\end{pro}

\noindent {\bf Proof.} a) Let $p,q\in
\Min(C)$ and $m\in \Spec(C)$ such that $p+q\subseteq m$. 
Since $C\subseteq
D$ satisfies LO and GD, there exist $P,Q,M\in \Spec(D)$ with 
$P\cap
C=p,Q\cap C=q,M\cap C=m$, and $P+Q\subseteq M.$ Choose 
$P_0,Q_0\in \Min(D)$
such that $P_0\subseteq P$ and $Q_0\subseteq Q$. Therefore $
P_0+Q_0\subseteq M$, with $P_0\cap A=p$ and $Q_0\cap A=q$. 
Since $D$
satisfies MPC, we have $P_0=Q_0$; consequently $p=q$, as 
desired.\\
b) It suffices to treat $A$. Now, $A\otimes_kB$ is $A$-
flat, and so
$A\rightarrow $$A\otimes_kB$
satisfies GD. It also satisfies LO by [26, Proposition 
2.3]. Apply
(a), to complete the proof. $\diamondsuit$\bigskip

The following example shows that Vamos' result (mentioned 
above) does not
extend to arbitrary $k$-algebras. It also provides a 
counterexample
to the converse of Proposition 3.1(b).
\begin{example} There exist a separable 
algebraic field
extension
$K$ of finite degree over $k$ and a $k$-algebra $A$ 
satisfying MPC such that
$K\otimes_kA$ fails to satisfy MPC.\end{example}

Let $k=I\!\!R$ and $K={\bf C}$ be the fields of real numbers
and complex numbers, respectively. Let
$V:={\bf C}[X]_{(X)}={\bf C}+X{\bf C}[X]_{(X)}$ and $A:= I\!
\!R + X{\bf
C}[X]_{(X)}$.
Clearly, $A$ is a one-dimensional local domain with 
quotient field
$L={\bf C}(X)$ and maximal ideal $p= X{\bf C}[X]_{(X)}$, 
such that $\frac
Ap=I\!\!R$.
We wish to show that $K\otimes_{I\!\!R}A$ does not satisfy 
MPC.
Indeed, let $f(Z)=Z^2+1$ be the minimal polynomial of $i$ 
over $I\!\!R$. We
have $K\otimes_{I\!\!R}A\cong \frac
{I\!\!R[Z]}{(f(Z))}\otimes_{I\!\!R}A\cong \frac {A[Z]}{(f
(Z))}$.
Therefore, the minimal prime ideals of $K\otimes_{I\!\!R}A$ 
are
$\overline I=\frac I{(f)}$ and
$\overline J=\frac J{(f)}$, where
$I=(Z-i)L[Z]\cap A[Z]$ and $J=(Z+i)L[Z]\cap A[Z]$. Since
$K\otimes_{I\!\!R}A$
is an integral extension of $A$, then so are
$\frac {A[Z]}I\cong \frac {A[Z]/(f)}{\overline I}$ and
$\frac {A[Z]}J\cong \frac {A[Z]/(f)}{\overline J}$, whence
$\dim(\frac {A[Z]}I) = \dim(\frac {A[Z]}J)= \dim (A) =1$. 
It follows that
$I$
and $J$ are not maximal ideals in $A[Z]$. Then, there exist 
$P_I$ and $P_J$
in
$\Spec(A[Z])$ such that $I\subset P_I$ and $J\subset P_J$. 
Clearly, $P_I\cap
A=P_J\cap A=p$.
Further, since $f\in I \cap J$ and $f \not\in p[Z]$, then
$P_I$ and $P_J$ are both uppers to $p$. As $\frac Ap=I\!\!R$
and $f$ is an irreducible monic polynomial over $I\!\!R$, 
it follows
that $P_I=P_J=(p,f)$ (cf. [18, Theorem 28]). Therefore
$I+J\subseteq P:=(p,f)$, and hence $\overline I+\overline 
J\subseteq
\overline P:=\frac P{(f(Z))}$.
Consequently, $\frac {A[Z]}{(f)}\cong K\otimes_{I\!\!R}A$ 
does not
satisfy MPC, establishing the claim. $\diamondsuit$\bigskip

We next investigate various contexts for the tensor product 
to
inherit the MPC property. The following result treats the 
case where
the ground field $k$ is algebraically closed.
\begin{thm} Let $k$ be an algebraically 
closed field. Let
$A$ and
$B$ be $k$-algebras. Then $A\otimes _kB$ satisfies MPC if 
and only if $A$
and $B$ each satisfy MPC.\end{thm}

\noindent {\bf Proof.} Proposition 3.1(b) handles the 
``only if" assertion.
Next,
assume that $A$ and $B$ each satisfy MPC. Let $P_0,Q_0\in
\Min(A\otimes_kB)$ and $P\in \Spec(A\otimes _kB$ $)$ such 
that $%
P_0+Q_0\subseteq P$.
Let $p_1:=P_0\cap A,q_1:=P_0\cap B$ and $p_2:=Q_0\cap
A,q_2:=Q_0\cap B$. We have $p_1,p_2\in \Min(A)$ and $q_1,q_2
\in \Min(B)$,
since $A\subseteq A\otimes_kB$ and $B\subseteq A\otimes_kB$ 
each satisfy GD.
Let $p:=P\cap A$ and $q:=P\cap B$. Then $p_1+p_2\subseteq 
p$ and $
q_1+q_2\subseteq q$. As $A$ and $B$ each satisfy MPC, 
$p_1=p_2=:p_0$ and
$q_1=q_2=:q_0$. Since $k$ is algebraically closed, it 
follows from
[27, Corollary 1, Ch. III, p. 198] and the lattice-
isomorphism in
[26, Proposition 2.3] that there is a unique prime $Q$ of 
$A\otimes_kB$ that
is minimal with respect to the properties $Q\cap A=p_0$, 
$Q\cap B=q_0$.
Hence, $P_0=Q=Q_0$, and the proof is complete. 
$\diamondsuit$\bigskip

The next theorem generalizes the above-mentioned result of 
Vamos.
\begin{thm} If $A$ and $B$ are integrally 
closed domains
that
are $k$-algebras, then $A\otimes _kB$ satisfies MPC.\end{thm}

\noindent {\bf Proof.} Let $K$ (resp., $L$) denote the 
quotient field of $A$
(resp., $B$). Let $K_s$ (resp., $L_s$) denote the separable 
algebraic
closure of $k$ in $K$
(resp., in $L$). Since $A$ is integrally closed and 
$k\subseteq A\subseteq
K$,
the algebraic closure of $k$ in $K$ is contained in $A$. In 
particular,
$K_s\subseteq A$; and, similarly, $L_s\subseteq B$. By [25, 
Theorem 3],
$\Min(K\otimes _kL)$ and $%
\Spec(K_s\otimes _kL_s)$ are canonically homeomorphic, with 
the prime
ideals of $K_s\otimes
_kL_s$ being the contractions of the minimal prime ideals 
of $K\otimes_kL$.
Observe that $K\otimes_kL$ is the localization of 
$A\otimes_kB$ at
$\{a\otimes b:a\in A \setminus \{0\}, b\in B \setminus \{0\}
\}$. It follows
that
there is a one-to-one correspondence between $\Min(K\otimes 
_kL)$
and $\Min(A\otimes _kB)$. Since
$K_s\otimes _kL_s\subseteq A\otimes _kB\subseteq 
K\otimes_kL$, we
obtain, via contraction, a bijection between $\Min(A\otimes 
_kB)$
and $\Spec(K_s\otimes _kL_s)$. Now, consider $P_0,Q_0\in \Min
(A\otimes _kB)$
and $%
P\in \Spec(A\otimes _kB)$ such that $P_0+Q_0\subseteq P$. 
Taking contractions
to $K_s\otimes _kL_s$, we obtain $P_0^c=P^c$ and 
$Q_0^c=P^c$, since $\dim
(K_s\otimes _kL_s)=0$ [26]. In particular, $P_0^c=Q_0^c$. 
By the
above bijection, $P_0=Q_0$, as desired. $\diamondsuit$\bigskip

The proof of Theorem 3.4 actually gives the following 
result. Let
$A$ and $B$ be domains that are $k$-algebras. Let $K_s$ 
(resp., $L_s$)
be the separable algebraic closure of $k$ in the quotient 
field $K$
(resp., $L$) of $A$ (resp., $B$). If $K_s\subseteq A$ and 
$L_s\subseteq B$,
then $A\otimes _kB$ satisfies
MPC.\bigskip

Moving beyond the contexts of Theorems 3.3 and 3.4, we next 
show that
$A\otimes_kB$ can satisfy MPC when $k$ is not algebraically 
closed
and when $A$, $B$ are not integrally closed domains.
\begin{example} Let $k:=I\!\!\!Q$ be the field 
of rational
numbers and let $A:=B:=I\!\!\!Q(i)[X^2,X^3]$. The quotient 
field of $A$
(resp., $B$) is $K=L=I\!\!\!Q(i)(X)$. We can easily check 
that $A$ and
$B$ are not integrally
closed (in fact they are not seminormal), and $K_s=L_s=I\!\!
\!Q(i)$,
since $I\!\!\!Q(i)[X^2,X^3]\subseteq I\!\!\!Q(i)[X]$ which
is integrally closed. Then $K_s\subseteq A$ and 
$L_s\subseteq B$.
By the above remark, $A\otimes_{I\!\!\!Q}B$
satisfies MPC, although $k=I\!\!\!Q$ is not algebraically 
closed and $A,B$
are not
integrally closed. $\diamondsuit$\end{example}

In Example 3.2, we exhibited a separable algebraic extension
field $K$ of $k$ and a $k$-algebra $A$ satisfying MPC such 
that
$K\otimes_kA$ fails
to satisfy MPC. The following result studies the case where 
$K$ is
purely inseparable over $k$.
\begin{pro} Let $A$ be a $k$-algebra 
and $K$ a
purely inseparable field extension of $k$. Then $K\otimes 
_kA$ satisfies
MPC if and only if $A$ satisfies MPC.\end{pro}

\noindent {\bf Proof.} Proposition 3.1(b) handles the 
``only if"
assertion. Conversely, assume that $A$ satisfies MPC. Let 
$P_0,Q_0$
be minimal
prime ideals of $K\otimes _kA$ and let $P\in \Spec(K\otimes 
_kA)$ such
that $P_0+Q_0\subseteq P$. Put $p_0:=P_0\cap A,q_0:=Q_0\cap 
A$, and
$p:=P\cap A$.
Hence $p_0+q_0\subseteq p$. Of course, $p_0$ and $q_0$ are 
in $\Min(A)$ since
flatness
ensures that $A\subseteq K\otimes_kA$ satisfies GD. Thus, 
since $A$
satisfies MPC, we obtain
$p_0=q_0$. However,
$\Spec(K\otimes_kA)\rightarrow \Spec(A)$ is an injection, 
since
``radiciel" is a universal property [13]. Consequently,
$P_0=Q_0$, as desired. $\diamondsuit$\bigskip

Theorem 3.9 is an extension to tensor products of $k$-
algebras of the
result [11, Proposition 2.1] that $A[X_1,...,X_n]$
is an S-domain, for any domain $A$ and any integer $n\geq 1
$. This latter
result was generalized to infinite sets of indeterminates 
in [10, Corollary
2.13].\bigskip

First we establish the following preparatory lemmas.
\begin{lem} If $A$ is a ring that satisfies 
MPC, then $
A[X_1,...,X_n]$ is an S-ring, for every integer $n\geq 1$.\end{lem}

\noindent {\bf Proof.} By Remark 2.2(c), $A[X_1,...,X_n]$ 
satisfies MPC.
Thus, it suffices to show the result when $n=1$. Let $Q\in 
\Min(A[X])$. Then
there exists $q\in \Min(A)$ such that $Q=q[X]$. Hence,
$\frac{A[X]}Q\cong \frac Aq[X]$ is an S-domain, by the above
remark, since $\frac Aq$ is an integral domain. Thus, $A[X]$
is an S-ring. $\diamondsuit$
\begin{lem} Let $A$ be a $k$-algebra and let 
$K$ be a
field extension
of $k$ such that $K\otimes _kA$ satisfies MPC. Then 
$K\otimes _kA$ is an
S-ring if and only if either $A$ is an S-ring or $t.d.(K:k)
\geq 1.$\end{lem}

\noindent {\bf Proof.} Suppose that $t:=t.d.(K:k)=0$, 
i.e., that $K$ is
algebraic over $k$. Then $K\otimes _kA$ is an integral 
extension of $A$ and
thus satisfies LO. Furthermore $A\subseteq K\otimes _kA$ 
satisfies GD and so
it follows easily that $A$ inherits MPC from $K\otimes_kA$. 
It remains to
show
that if $P\in \Min(K\otimes_kA)$ and $p=P\cap A$, then 
$\frac {K\otimes
_kA}{P}$
is an S-domain if and only if $\frac Ap$ is an S-domain. 
The ``only if"
statement
follows from the proof of [19, Theorem 4.6], while the 
treatment of the
``if"
statement is similar to that of the proof of [19, Theorem 
4.9].

In the remaining case, $t:=t.d.(K:k)\geq
1$. Let $B$ be a transcendence basis of $K$ over $k$. As 
$K\otimes
_kA\cong K\otimes _{k(B)}(k(B)\otimes _kA)$, we see that $k
(B)\otimes _kA$
satisfies MPC, by Proposition 3.1(b). Also, if $X\in B$,
$B_1:=B\setminus\{X\}$ and
$S:=k(B_1)[X]\setminus\{0\}$, then $k(B)\otimes _kA
=k(B_1)(X)\otimes _kA\cong S^{-1}((k(B_1)\otimes _kA)[X])$. 
As
$k(B_1)\otimes _kA$ satisfies
MPC, Lemma 3.7 yields that $(k(B_1)\otimes_kA)[X]$ is an S-
ring. Hence, so
is
its ring of fractions $k(B)\otimes_kA$. Therefore, by the 
first case, so
is $K\otimes_{k(B)}(k(B)\otimes _kA)\cong K\otimes _kA$, to 
complete the
proof. $\diamondsuit$

\begin{thm} Let $A$ and $B$ be $k$-
algebras such that
$A\otimes
_kB$ satisfies MPC. Then $A\otimes _kB$ is an S-ring if and 
only if at least
one of
the following statements is satisfied:

1) $A$ and $B$ are S-rings;

2) $A$ is an S-ring and $t.d.(\frac Ap:k)\geq 1$ for each 
$p\in \Min(A)$;

3) $B$ is an S-ring and $t.d.(\frac Bq:k)\geq 1$ for each
$q\in \Min(B)$;

4) $t.d.(\frac Ap:k)\geq 1\;$and $t.d.(\frac Bq:k)\geq 1$ 
for each $p\in
\Min(A)$ and $q\in \Min(B)$.\end{thm}

\noindent {\bf Proof.} We claim that $A\otimes _kB$ is an S-
ring
if and only if $k_A(p)\otimes _kB$ and $A\otimes _kk_B(q)$ 
are S-rings for
each $p\in \Min(A)$ and $q\in \Min(B)$. Indeed,
assume
that $A\otimes _kB$ is an S-ring. Clearly, by [26, 
Proposition 2.3], for
each
minimal prime ideal $p$ of $A$, $\frac Ap\otimes _kB\cong 
\frac
{A\otimes_kB}{p\otimes_kB}$
satisfies MPC, and thus so does its ring of fractions $k_A
(p)\otimes _kB$.
Similarly, so does $A\otimes _kk_B(q)$, for each minimal 
prime ideal $q$ of
$B$.
In view of Remark 2.2(a), we may focus on $(P_2)$. Let 
$p\in \Min(A)$
and $P\in \Spec(A\otimes _kB)$ such that $P\cap A=p$ and $ht
(\frac P{p\otimes
_kB})=1$.
Since $p\in \Min(A),$ we have $\htt(P)=\htt(\frac P{p\otimes 
_kB})=1$. By
the hypothesis on $A\otimes_kB$, $1=
\htt(P[X])=\htt(\frac P{p\otimes _kB}[X])$. Hence, $k_A(p)
\otimes_kB$ is an
S-ring for each $p\in \Min(A)$. Similarly, so is
$A\otimes _kk_B(q)$ for each $q\in \Min(B)$.

Conversely, suppose that $k_A(p)\otimes _kB$ and $A\otimes 
_kk_B(q)$ are
S-rings for each $p\in \Min(A)$ and $q\in \Min(B)$. Let $P\in
\Spec(A\otimes _kB)$ such that $\htt(P)=1$. By [26, Corollary 
2.5], we have
that
either $p:=P\cap A$ is a minimal prime ideal of $A$ or 
$q:=P\cap B$ is a
minimal prime ideal of $B$. Without loss of generality,
$p\in \Min(A).$ Then $\htt(\frac
P{p\otimes _kB})=\htt(P)=1$. Since $k_A(p)\otimes _kB$ is an 
S-ring, we have
$1=
\htt(\frac P{p\otimes _kB}[X])=\htt(P[X])$. Consequently, 
$A\otimes _kB$ is an
S-ring, and the claim has been proved. The theorem now 
follows
from Lemma 3.8. $\diamondsuit$\bigskip

It is clear from the above proof that the statement of 
Theorem 3.9
remains true without the MPC hypothesis if we substitute 
$(P_2)$ for
the S-ring property.
\begin{cor} Let $k$ be an algebraically 
closed field.
Let $A$
and $B$ be domains that are $k$-algebras. Then $A\otimes 
_kB$ is an S-domain
if and
only if at least one of the following statements is 
satisfied:

1) $A$ and $B$ are S-domains;

2) $A$ is an S-domain and $t.d.(A:k)\geq 1$;

3) $B$ is an S-domain and $t.d.(B:k)\geq 1$;

4) $t.d.(A:k)\geq 1$ and $t.d.(B:k)\geq 1$.\end{cor}

\noindent {\bf Proof.} Apply Theorem 3.9, bearing in mind 
that
$A\otimes _kB$ is an integral
domain (hence satisfies MPC) since $k$ is algebraically 
closed
[27, Corollary 1, Ch. III, p. 198]. $\diamondsuit$
\begin{cor} Let $A$ and $B$ be 
integrally
closed domains that are $k$-algebras. Then $A\otimes _kB$ 
is an
S-ring if and only if at least one of the following
statements is satisfied:

1) $A$ and $B$ are S-domains;

2) $A$ is an S-domain and $t.d.(A:k)\geq 1$;

3) $B$ is an S-domain and $t.d.(B:k)\geq 1$;

4) $t.d.(A:k)\geq 1$ and $t.d.(B:k)\geq 1$.\end{cor}

\noindent {\bf Proof.} Combine Theorems 3.9 and 3.4. 
$\diamondsuit$
\end{section}\bigskip

\begin{section}{Strong S-property and Catenarity}

Our purpose in this section is to seek conditions for the 
tensor product
of two $k$-algebras to inherit the (stably) strong S-
property and
(universal) catenarity. The main theorem of this section 
generates new
families of
stably strong S-rings and universally catenarian rings. Our 
interest is
turned essentially to studying $A\otimes_kB$ in case at least one of 
$A$, $B$ is a
field extension
of $k$. Beyond this context, the study of these properties
becomes more intricate, as one may expect. In fact, a 
glance ahead to
Example 5.5 reveals a non-catenarian ring of the form 
$A\otimes_kB$ in
which $A$, $B$ are each universally catenarian domains (in 
fact DVRs).\bigskip

To determine when a tensor product of $k$-algebras is 
catenarian, we first
need
to know when it is LFD. That is handled by the first result 
of this
section.
\begin{pro}  Let $A$ and $B$ be k-algebras. Then:\\
a) If $A\otimes_kB$ is LFD, then so are $A$ and $B$, and  either $t.d.(\frac Ap:k) < \infty$ for each prime ideal $p$ of $A$ or $t.d.(\frac Bq:k) < \infty$ for each prime ideal $q$ of $B$.\\
b) If both $A$ and $B$ are LFD and either $t.d.(A:k) < \infty$ or $t.d.(B:k) < \infty$, then $A\otimes_kB$ is LFD. The converse holds provided $A$ and $B$ are domains. \end{pro}

The proof of this proposition requires the following 
preparatory lemma.
\begin{lem} Let $K$ and $L$ be field
extensions of $k$. Then $K\otimes_kL$ is LFD if and only if 
either $
t.d.(K:k)< \infty $ or $t.d.( L:k)< \infty$.\end{lem}

\noindent {\bf Proof.} $\Leftarrow )$ Straightforward,
since $\dim(K\otimes_kL)=\min(t.d.( K:k),t.d.( L:k))$ (cf. 
[23, Theorem
3.1]).

$\Rightarrow )$ Let $B$ (resp., $B^{\prime }$) be a 
transcendence basis
of $K$ (resp., $L$) over $k$. As
$K\otimes_kL\cong
K\otimes_{k(B)}(k(B)\otimes_kk(B^{\prime}))\otimes_{k(B^
{\prime})}L$,
then $k( B)\otimes_kk(B^{\prime})\subset_{\rightarrow}
K\otimes_kL$ is an integral extension that satisfies
GD. Therefore, $K\otimes_kL$ is LFD if and only if
$k(B)\otimes_kk(B^{\prime})$ is LFD. Suppose that $t.d.(K:k)
=t.d.(L:k)
=\infty$. Let $T:=k(x_1, x_2, ...)
\otimes_kk(y_1, y_2, ...)$, where the $x_i\in B$ and the 
$y_i\in
B^{\prime}$. Since $T\subseteq k(B)\otimes_{k(x_1, 
x_2, ...)}T$ and
$k(B)\otimes_{k(x_1, x_2, ...)}T\subseteq (k(B)\otimes_{k
(x_1, x_2, ...)}T)
\otimes_{k(y_1, y_2, ...)}k(B')\cong k(B)\otimes_kk(B^
{\prime})$
are ring extensions that satisfy GD and LO, then so does 
$T\subseteq
k(B)\otimes_kk(B^{\prime})$.
Thus $T$ is not LFD $\Rightarrow
k(B)\otimes_kk(B^{\prime })$ is not LFD $\Rightarrow 
K\otimes_kL$
is not LFD.

Let $K_{n}=k(x_1,...,x_n)$ and
$S_n = k[y_1,...,y_n]\setminus\{ 
0\}$,
for each $n\geq 1$. Consider the following ring 
homomorphisms:
\[K_{n}[y_1,...,y_n]\subset \stackrel{i_n}{_
{\longrightarrow }}
k(x_1, x_2, ...)[y_1, y_2, ...]
\stackrel{\varphi }{\longrightarrow } k(x_1, x_2, ...)\]
where $\varphi (y_i)=x_i$ for $i\geq 1$. Let $M=Ker(\varphi)
$ and $%
M_n=M\cap K_{n}[y_1,...,y_n]=Ker(\varphi _n)$,
where $\varphi _n:=\varphi \circ i_n$, for all $n\geq 1$. 
Since $x_1,...,x_n$
are algebraically independent over $k$, $M_n\cap S_{n}
=\emptyset$, for all
$n\geq 1$. On the other hand, since $K_{n}[y_1,...,y_n]$ is 
an AF-domain (we recall early in 
Section 5 the definition of an AF-domain), then, for 
every $n\geq 1$,
\[\htt(M_n)+t.d.(\frac{K_{n}[y_1,...,y_n]}{M_n}
:K_{n})=t.d.(K_{n}[y_1,...,y_n]:K_{n})=n.\]
Hence $\htt(M_n)=n$, since $\frac{K_{n}[y_1,...,y_n]}{M_n}%
\cong K_{n}$, for all $n\geq 1$. Therefore $M\cap 
S=\emptyset$, where
$S:=\bigcup_{n}S_n=k[y_1, y_2, ...]\setminus\{0\}$. We wish 
to show that 
$\htt(M)=\infty$.
Indeed, observe that, for any integer $n\geq 1$,
\begin{eqnarray*}
M_nk(x_1, ...)[y_1, ...] & = & M_n(k(x_1,
...)\otimes_{K_{n}}K_{n}[y_1,...,y_n]\otimes_kk[y_
{n+1},...])\\
& = & k(x_1,...)\otimes_{K_{n}}M_n\otimes_kk[y_
{n+1},...], \
\hbox{and}\\
\frac {k(x_1, ...)[y_1, ...]}{M_nk(x_1, ...)
[y_1, ...]}
& \cong & \frac {k(x_1,...)\otimes_{K_{n}}K_{n}
[y_1,...,y_n]
\otimes
_kk[y_{n+1},...]}{k(x_1, ...)\otimes_{K_{n}}M_n
\otimes_kk[y_{n+1},...]}\\
&\cong & k(x_1,...)\otimes_{K_{n}}\frac {K_{n}
[y_1,...,y_n]}{M_n}
\otimes_kk[y_{n+1},...]\ \hbox{(cf. [26])}\\
&\cong & k(x_1,...)\otimes_{K_{n}}K_{n}\otimes_kk[y_
{n+1},...] \\
&\cong & k(x_1, ...)[y_{n+1},...],\
\hbox{an integral domain.}
\end{eqnarray*}

\noindent Thus $M_nk(x_1, ...)[y_1, ...]$ is
a prime ideal in $k(x_1, ...)[y_1, ...]$, for all
$n\geq 1$. Since $K_{n}[y_1,...,y_n]\rightarrow
k(x_1, ...)[y_1, ...]$ is a faithfully flat 
homomorphism (and
hence satisfies GD),
we obtain $M_nk(x_1, ...)[y_1, ...]
\cap K_{n}[y_1,...,y_n]=M_n$, and thus
$\htt(M_nk(x_1, ...)[y_1, ...])\geq \htt(M_n)=n$. By 
direct limits (cf.
[10]),
it follows that $\htt(M)=\infty$, as desired. Consequently, 
$S^{-1}M$
is a prime ideal of $T = k(x_1, ...)\otimes _kk(y_1, ...)$ with
$\htt(S^{-1}M)=\infty$. Therefore $T$ is not LFD, completing 
the
proof.$\diamondsuit$\bigskip

\noindent {\bf Proof of Proposition 4.1.} a) Assume that
$A\otimes _kB $ is LFD. Let $p\in \Spec(A)$ 
and $q\in \Spec(B)$. As the extensions $A\subseteq 
A\otimes_kB$ and $B\subseteq A\otimes_kB$ satisfy LO, 
there exist prime ideals $P$ and $Q$ of $A\otimes_kB$ such 
that $P\cap A=p$ and $Q\cap B=q$. By [26, Corollary 2.5], 
$\htt(p)\leq \htt(P)<\infty$ and $\htt(q)\leq \htt(Q)<\infty$. It 
follows that $A$ and $B$ are LFD. Now, suppose that there exists a prime 
ideal $q$ of $B$ such that $t.d.(\frac Bq:k)=\infty$. Let 
$p$ be any prime ideal of $A$. Then $\frac 
Ap\otimes_k\frac Bq\cong \frac {A\otimes_kB}
{p\otimes_kB+A\otimes_kq}$ is LFD. Hence $k_A(p)\otimes_kk_B
(q)$ is LFD, since it is a ring of fractions of $\frac 
Ap\otimes_k\frac Bq$. Therefore, by Lemma 4.2, $t.d.(k_A
(p):k)=t.d.(\frac Ap:k)<\infty$.  

\noindent b) Suppose that $t.d.(A:k)<\infty$ and both 
$A$ and $B$ are
LFD. Consider a chain $\Omega:=\{P_0\subset P_1
\subset ...\subset P\}$ of
prime ideals of $A\otimes _kB$
and let $l$ be its length. We claim that $l$ is finite, with an upper
bound depending on $P$.
Let $p_0\subset ...\subset p_r=p:=P\cap A$ and $q_0
\subset ...\subset
q_s=q:=P\cap B$
be the chains of intersections of $\Omega $ over $A$ and 
$B$, respectively.
We can
partition $\Omega
$ into subchains $\Omega _{ij}$ the prime ideals of which 
contract to $p_i$
in $\Spec(A)$ and $q_j$ in $\Spec(B)$. Thus each $\Omega _{ij}
$ is of length
$l_{ij} \leq
\dim (k_A(p_i)\otimes _kk_B(q_j))$, by [26, Proposition 
2.3].
Therefore, we have
\begin{eqnarray*}
l &\leq & \sum_{i=0,j=0}^{r,s}(\dim (k_A(p_i)
\otimes _kk_B(q_j))+1)
\\
 &\leq & \sum_{i=0,j=0}^{r,s}(min(t.d.(\frac A
{p_i}:k),t.d.(\frac
B{q_j}:k))+1)\ \hbox{[23, Theorem 3.1]}\\
 &\leq & \sum_{i=0,j=0}^{r,s}(t.d.(\frac A
{p_i}:k)+1)\\
&\leq & (t.d.(A:k)+1)(r+1)(s+1)\\
&\leq & (t.d.(A:k)+1)(\htt(p)+1)(\htt(q)+1) <\infty,\ \hbox{as desired}.
\end{eqnarray*}

Now, if $A$ and $B$ are domains, then the converse holds, by (a). $\diamondsuit$\bigskip

Given an integer $n\geq 1$, Malik and Mott proved that $A
[X_1,...,X_n]$ is
a strong S-ring if and only if so is $A_p[X_1,...,X_n]$
for each prime ideal $p$ of $A$ [19, Theorem 3.2]. We next
extend this result to tensor products of $k$-algebras.
\begin{pro} Let $A_1$ and $A_2$ be $k$-
algebras.
Then the
following statements are equivalent:

1) $A_1\otimes _kA_2$ is a strong S-ring (resp., 
catenarian);

2) $S_1^{-1}A_1\otimes _kS_2^{-1}A_2$ is a strong S-ring
(resp., catenarian)
for each multiplicative subset $S_i$ of $A_i,$ for $i=1,2;$

3) $(A_1)_{p_1}\otimes _kA_2$ is a strong S-ring
(resp., catenarian) for
each $p_1\in \Spec(A_1)$;

4) $(A_1)_{m_1}\otimes _kA_2$ is a strong S-ring (resp., 
catenarian) for
each $m_1\in \Max(A_1)$;

5) $A_1\otimes _k(A_2)_{p_2}$ is a strong S-ring (resp., 
catenarian) for
each $p_2\in \Spec(A_2)$;

6) $A_1\otimes _k(A_2)_{m_2}$ is a strong S-ring (resp., 
catenarian) for
each $m_2\in \Max(A_2)$;

7)$(A_1)_{m_1}\otimes _k(A_2)_{m_2}$ is a strong S-ring 
(resp.,
catenarian) for each $m_i\in \Max(A_i)$, for $i=1,2$.\end{pro}

\noindent {\bf Proof.} The class of strong S- (resp., 
catenarian) rings
is stable under formation of rings of fractions. Thus
$(1)\Rightarrow (2)\Rightarrow (3)\Rightarrow (4)
\Rightarrow (7)$,
and $(2)\Rightarrow (5)\Rightarrow (6)\Rightarrow (7)$. 
Therefore,
it suffices to prove that $(7)\Rightarrow (1)$. Note that
if $(A_1)_{m_1}\otimes _k(A_2)_{m_2}$
satisfies MPC for each maximal ideal $m_i$ of $A_i$, for 
$i=1,2$,
then $A_1\otimes _kA_2$ satisfies MPC. Indeed, let
$P_1$ and $P_2$ be two minimal prime ideals contained in a 
common prime
ideal $P$ of $A_1\otimes _kA_2$. Choose a maximal ideal 
$m_i$ of $A_i$ such
that $P\cap A_i\subseteq m_i$, for $i=1,2$. Then $P_i((A_1)_
{m_1}\otimes
_k(A_2)_{m_2})\subseteq P((A_1)_{m_1}\otimes _k(A_2)_{m_2})
$, for $i=1,2$.
Hence, by hypothesis, $P_1((A_1)_{m_1}\otimes
_k(A_2)_{m_2})=P_2((A_1)_{m_1}\otimes _k(A_2)_{m_2})$. 
Taking contractions
to $A_1\otimes _kA_2$, we obtain $P_1=P_2$, since $(A_1)_
{m_1}\otimes
_k(A_2)_{m_2}$ is a ring of fractions of $A_1\otimes _kA_2
$. Then
$A_1\otimes_kA_2$
satisfies MPC. Also, if $(A_1)_{m_1}\otimes _k(A_2)_{m_2}$ 
is LFD for
each maximal ideal $m_i$ of $A_i$, for $i=1,2$, then it is 
clear that
$A_1\otimes_kA_2$ is LFD.

Now suppose that (7) holds. Let $P\subset Q$ be a saturated 
chain in
$\Spec(A_1\otimes _kA_2)$, $p_i:=P\cap A_i$ and $q_i:=Q\cap 
A_i$, for
$i=1,2$. Choose $
m_i\in \Max(A_i)$ such that $p_i\subseteq q_i\subseteq m_i$, 
for $i=1,2$.
Then
$P((A_1)_{m_1}\otimes _k(A_2)_{m_2})\subset Q((A_1)_{m_1}
\otimes
_k(A_2)_{m_2})$ is a saturated chain in $\Spec((A_1)_{m_1}
\otimes
_k(A_2)_{m_2})$.
Since $(A_1)_{m_1}\otimes _k(A_2)_{m_2}$ is a strong S-ring 
(resp.,
catenarian), we have $P((A_1)_{m_1}\otimes _k(A_2)_{m_2})[X]
\subset
Q((A_1)_{m_1}\otimes _k(A_2)_{m_2})[X]$ is a saturated 
chain in $\Spec(
(A_1)_{m_1}\otimes _k(A_2)_{m_2}[X])$ (resp., $\htt(Q((A_1)_
{m_1}\otimes
_k(A_2)_{m_2}))=1+\htt(P((A_1)_{m_1}\otimes _k(A_2)_{m_2}))
$. Therefore
$ht\frac{Q[X]}{P[X]}=1$ (resp., $\htt(Q)=1+\htt(P)$). Then (1) 
holds, completing
the proof. $\diamondsuit$\bigskip

It will follow from Theorem 4.9 (proved below) that if $K$ 
and $L$ are field
extensions
of $k$ with $t.d.(K)<\infty$, then $K\otimes_kL$ is a 
strong S-ring and
catenarian.
Applying Proposition 4.3, it follows that if $A$ and $B$
are von Neumann regular $k$-algebras with $t.d.(A)<\infty$, 
then
$A\otimes_kB$ is a strong S-ring and catenarian.
\begin{pro} Let $A$ be a $k$-algebra 
and $K$
an algebraic field
extension of $k$. If $K\otimes _kA$ is a strong S-ring 
(resp.,
catenarian), then $A$ is a strong S-ring (resp., 
catenarian).\end{pro}

\noindent {\bf Proof.} The strong S-property is 
straightforward from
Proposition 2.3. Assume that $K\otimes _kA$ is catenarian. 
Then $K\otimes
_kA$ satisfies MPC, and thus, by Proposition 3.1(b), $A$ 
satisfies MPC. Let
$p\subset q$
be a saturated chain of prime ideals of $A$. Since 
$K\otimes _kA$ is an
integral extension of $A$, there exists a saturated chain 
of prime ideals $%
P\subset Q$ of $K\otimes _kA$ such that $P\cap A=p$ and 
$Q\cap A=q$. Hence
$\htt(Q)=1+\htt(P)$. As $A\subset _{\rightarrow }K\otimes
_kA$ satisfies also GD, we obtain $\htt(q)=\htt(Q)=1+\htt(P)=1+ht
(p)$. Since, by
Proposition 4.1, $A$ is LFD, we conclude that $A$ is 
catenarian.
$\diamondsuit$\bigskip

Note that Proposition 4.4 fails, in general, when the 
extension
field $K$ is no longer algebraic over $k$, as it is shown 
by Example 5.2 and
Example 5.3.\bigskip

Next, we investigate sufficient conditions, on a $k$-
algebra $A$ and a field
extension $K$ of $k$, for $K\otimes_kA$ to inherit the 
(stably) strong
S-property
and (universal) catenarity.
\begin{pro} Let $A$ be a $k$-algebra 
and $K$ a purely
inseparable field extension of $k$. Then $K\otimes _kA$ is
a strong S-ring (resp., stably
strong S-ring, catenarian, universally catenarian) if and 
only if so is
$A$.\end{pro}

\noindent {\bf Proof.} $k\subset_{\rightarrow} K$ is 
radiciel,
hence a universal homeomorphism. In particular, both
$A\subset_{\rightarrow} K\otimes_kA$ and (for each $n\geq 1
$)
$A[X_1,...,X_n]
\subset_{\rightarrow} K\otimes_kA[X_1,...,X_n]\cong
(K\otimes_kA)[X_1,...,X_n]$
induce order-isomorphisms on Specs. Moreover, by 
Proposition 3.6,
$K\otimes_kA$
satisfies MPC if and only if $A$ satisfies MPC. Hence, the 
``catenarian" and
``universally catenarian" assertions now follow 
immediately. Also, by
applying Spec to the commutative diagram
$$
\begin{array}{ccc}
A & \subset_{\longrightarrow} & K\otimes_kA\\
\downarrow & & \downarrow\\
A[X] & \subset_{\longrightarrow}& K\otimes_kA[X]
\end{array}
$$
\noindent we obtain the ``strong S-ring" assertion and, 
hence,
the ``stably strong S-ring" assertion. $\diamondsuit$
\begin{pro} Let $A$ be a domain that 
is a $k$-algebra
and $K$ an
algebraic field extension of $k$. Assume that $A$ contains 
a separable
algebraic
closure of $k$. Then $K\otimes _kA$ is a strong S-ring 
(resp., stably strong
S-ring, catenarian, universally catenarian) if and only if 
so is $A$.\end{pro}

\noindent {\bf Proof.} Proposition 4.4 handles the ``only 
if" assertion.
Conversely, let $\overline k$ be the
separable algebraic closure of $k$ contained in $A$.
First, we claim that the contractions of any adjacent prime 
ideals
of $K\otimes_kA[X_1,...,X_n]$ are adjacent in $A
[X_1,...,X_n]$. Indeed, let
$n$
be a positive integer and $P\subset Q$ be a pair of 
adjacent prime ideals of
$K\otimes_kA[X_1,...,X_n]$. Put $P^{\prime}:=P\cap A
[X_1,...,X_n]$ and
$Q^{\prime}
:=Q\cap A[X_1,...,X_n]$. Since also $\overline k\subseteq
\frac {A[X_1,...,X_n]}{P^{\prime}}$, then $K\otimes_k\frac
{A[X_1,...,X_n]}{P^{\prime}}$
satisfies MPC (see the remark following Theorem 3.4). 
Furthermore, since $K$
is algebraic
over $k$, $P$ is the unique prime ideal of $K\otimes_kA
[X_1,...,X_n]$
contained in $Q$ and contracting to $P^{\prime}$ by [26, 
Proposition 2.3]
and [23, Theorem 3.1]. Hence,
$1=\htt(\frac QP)=\htt(\frac Q{K\otimes_kP^{\prime}})
= \htt(\frac {Q^{\prime}}{P^{\prime}})$, proving the claim. 
Now the ``strong
S-ring"
and ``stably strong S-ring" assertions follow easily. 
Moreover, since
$K\otimes_kA[X_1,...,X_n]$ is an integral extension of $A
[X_1,...,X_n]$
that satisfies GD, for any integer $n$, we have for any 
prime ideals
$P\subseteq Q$
of $K\otimes_kA[X_1,...,X_n]$, $\htt(P)=\htt(P^{\prime})$ and
$\htt(Q)=\htt(Q^{\prime})$,
where $P^{\prime}:=P\cap A[X_1,...,X_n]$ and $Q^
{\prime}:=Q\cap
A[X_1,...,X_n]$.
Then, in view of the above claim, the ``catenarian" and 
``universally
catenarian"
statements follow, completing the proof. $\diamondsuit$
\begin{thm} Let $A$ be a domain that is a 
$k$-algebra and
$K$ an algebraic
field extension of $k$. Assume that the
integral closure $A^{\prime}$ of $A$ is a Pr\"ufer domain. 
Then
$K\otimes_kA$ is a stably strong S-ring.\end{thm}

\noindent {\bf Proof.} We claim that $K\otimes_kA^{\prime}$ 
is a stably
strong S-ring. In fact, let $P_0$ be a minimal prime ideal 
of
$K\otimes_kA^{\prime}$.
Then $P_0\cap A^{\prime}=(0)$, and thus $\frac {K\otimes_kA^
{\prime}}{P_0}$
is an integral extension of $A^{\prime}$. Since $A^{\prime}
$ is a
Pr\"ufer domain,
$\frac {K\otimes_kA^{\prime}}{P_0}$ is a stably strong S-
domain by [19,
Proposition 4.18].
It follows that $K\otimes_kA^{\prime}$ is a stably strong S-
ring,
as desired. Proposition 2.3 completes the proof. 
$\diamondsuit$
\begin{thm} Let $A$ be an LFD
Pr\"ufer domain that is a $k$-algebra and $K$ an algebraic 
field extension
of $k$. Then $K\otimes _kA$ is catenarian.\end{thm}

\noindent {\bf Proof.} First, we have that $K\otimes_k\frac 
Ap$ satisfies
MPC, by Theorem
3.4, since $
\frac Ap$ is integrally closed for any $p\in \Spec(A)$. An 
argument similar
to the treatment of the claim in the proof of
Proposition 4.6 allows us to see that the contractions of 
any adjacent prime
ideals of $K\otimes_kA$ are adjacent in $\Spec(A)$. Then, 
since $K\otimes_kA$
is
an integral extension of $A$ that satisfies GD, the result 
follows, since
the contraction map from $\Spec(K\otimes_kA)$ to $\Spec(A)$
preserves height. $\diamondsuit$
\begin{thm} Let $A$ be a Noetherian domain 
that is a
$k$-algebra and $K$ a field extension of $k$ such that $t.d.
(K:k)<\infty$.
Then $K\otimes_kA$ is a stably strong S-ring.
\noindent If, in addition, $K\otimes _kA$ satisfies MPC and 
$A[X]$ is
catenarian, then $K\otimes _kA$ is universally catenarian.\end{thm}

\noindent {\bf Proof.} Recall first that a Noetherian ring $A$ is universally catenarian if and only if $A[X]$ is catenarian [22]. We have $K\otimes _kA\cong K\otimes
_{k(X_1,...,X_t)}S^{-1}A[X_1,...,X_t]$, where $t:=t.d.(K:k)
$ and $%
S:=k[X_1,...,X_t]\setminus\{0\}$. Since $S^{-1}A
[X_1,...,X_t]$ is
Noetherian,
it
suffices to handle the case where $K$ is algebraic over 
$k$. Thus, in
that case, $K\otimes _kA$ is an integral extension of a 
Noetherian domain
$A$. Let $P_0$ be a minimal prime ideal of $K\otimes_kA$. 
By GD, $P_0\cap
A=(0)$.
It follows that $\frac {K\otimes_kA}{P_0}$ is an integral 
extension of $A$.
Hence, by
[19, Proposition 4.20], $\frac {K\otimes_kA}{P_0}$ is a 
stably strong
S-domain,
whence $K\otimes_kA$ is a stably strong S-ring.
Now, assume that $K\otimes _kA$ satisfies MPC and $A[X]$ is 
catenarian. Let
$P_0$ be a minimal prime ideal of $K\otimes _kA$. As above, 
$\frac {K\otimes
_kA}{P_0}$
is an integral extension of $A$.
By [22, Theorem 3.8], $\frac {K\otimes _kA}{P_0}$ is a 
universally
catenarian
domain. It follows that $K\otimes _kA$ is a universally 
catenarian ring. The
proof is complete. $\diamondsuit$
\begin{cor} Let $K$ and $L$ be field 
extensions of $k$
such that $t.d.(K:k)< \infty $. Then $K\otimes _kL$ is 
universally
catenarian.\end{cor}

\noindent {\bf Proof.} $K\otimes _kL$ is LFD by Lemma 4.2, 
and satisfies MPC
by [25, Corollary 4]. Theorem 4.9 completes the proof. 
$\diamondsuit$
\begin{cor} Let $A$ be a one-
dimensional $k$-algebra and
$K$ an
algebraic field extension of $k$. Then the following 
conditions are
equivalent:

i) $K\otimes _kA$ is a stably strong S-ring;

ii) $K\otimes _kA$ is a strong S-ring;

iii) $A$ is a strong S-ring;

iv) $A$ satisfies $(P_2)$.

\noindent If, in addition, $K\otimes _kA$ satisfies MPC, 
then the following
conditions are equivalent and the assertions (i)-(iv)
are equivalent to each of (v)-(viii):

v) $K\otimes _kA$ is universally catenarian;

vi) $(K\otimes _kA)[X]$ is catenarian;

vii) $K\otimes _kA$ is an S-ring;

viii) $A$ is an S-ring.\end{cor}

\noindent {\bf Proof.} By Proposition 2.3 and Proposition 
2.5, we have
$(i)\Leftrightarrow (ii)\Rightarrow (iii)\Leftrightarrow 
(iv)$. \\
\noindent $iv)\Rightarrow ii)$ Assume that (iv) holds. Let 
$P$ be a minimal
prime ideal of $K\otimes_kA$ and $p:=P\cap A$. If $\frac 
{K\otimes_kA}P$
is a field, then it is an S-domain. If $\dim(\frac 
{K\otimes_kA}P)=1$,
then $\dim(\frac Ap)=1$ (since $\frac {K\otimes_kA}P$ is an 
integral
extension of $\frac Ap$); therefore $\frac Ap$ is an S-
domain by
Proposition 2.5, whence
$\frac {K\otimes_kA}P$ is an S-domain by [19, Theorem 4.2].
We conclude that $ K\otimes_kA$ is a strong S-ring. Thus 
the statements
(i)-(iv)
are equivalent.
On the other hand, the
assertions (v)-(vii) are equivalent, by Proposition 2.5. 
Also,
$(vii)\Leftrightarrow (viii)$, by Lemma 3.8. Apply 
Proposition 2.5 to
complete the proof. $\diamondsuit$
\begin{pro} Let $A$ be a two-
dimensional
$k$-algebra and $K$
an algebraic field extension of $k$ such that $K\otimes 
_kA$ satisfies MPC.
Then $K\otimes _kA$ is a strong S-ring (resp., catenarian) 
if and only if
so is
$A$.\end{pro}

\noindent {\bf Proof.} The ``only if" assertion follows 
from Proposition
2.3.
Conversely, we first show that the contractions of any pair 
of adjacent
prime ideals of $K\otimes_kA$ are adjacent in $\Spec(A)$. In 
fact, let
$P\subset Q$ be a pair of adjacent prime ideals in 
$K\otimes _kA$,
$p:=P\cap A$ and $q:=Q\cap A$. If $\htt(P)=1$, then $\htt(p)=1$ 
and hence
$\htt(\frac qp)=1$, since $\dim (A)=2$. In the remaining 
case, $P$ is a
minimal prime ideal of $K\otimes _kA$. Since $K\otimes _kA$
satisfies MPC, $P$ is the unique minimal prime ideal
contained
in $Q$. Then $\htt(Q)=\htt(\frac QP)=1$. It follows that
$\htt(\frac qp)\leq \htt(q)=\htt(Q)=1$, since $K\otimes_kA$ is an
integral extension of $A$ that satisfies GD. Then $\htt(\frac 
qp)=1$.
Hence, the ``strong S-ring" assertion follows immediately. 
As the
contraction map from $\Spec(K\otimes_kA)$ to $\Spec(A)$ 
preserves height,
the ``catenarian" assertion also holds. $\diamondsuit$

Next, we state the main theorem of this section. It 
generates
new families of stably strong S-rings and universally
catenarian rings.
\begin{thm} Let $A$ be an LFD $k$-algebra 
and $K$
a field extension
of $k$ such that either $t.d.(A:k)< \infty $ or
$t.d.(K:k)< \infty
$. Let $B$ be a transcendence basis of $K$ over $k$, and 
let $L$
be the
separable algebraic closure of $k(B)$ in $K$. Assume that
$[L:k(B)]<
\infty$. If $A$ is a stably strong S-ring (resp.,
universally
catenarian and $K\otimes _kA$ satisfies MPC), then 
$K\otimes _kA$
is a stably
strong S-ring (resp., universally catenarian).\end{thm}

The proof of this theorem requires the following preparatory
result.
\begin{pro} Let $A$ be an LFD $k$-
algebra and $K$ a
purely
transcendental field extension of $k$ such that either $t.d.
(A:k)< \infty $
or $%
t.d.(K:k)< \infty $. If $A$ is a stably strong S-ring 
(resp.,
universally catenarian), then $K\otimes _kA$ is a stably 
strong S-ring
(resp.,
universally catenarian).\end{pro}

\noindent {\bf Proof.} First note that the stably strong S-
property and
universal catenarity are stable under formation of
rings of fractions. Let $K=k(B)$,
where $B$ is a transcendence basis of $K$
over $k$. If $B$ is a finite set $\{X_1,...,X_n\}$, then
$K\otimes _kA\cong
S^{-1}A[X_1,...,X_n]$, where $S:=k[X_1,...,X_n]\setminus\{0
\}$. Clearly,
$K\otimes _kA$
is a stably strong S-ring (resp., universally catenarian), 
if
$A$ is. Hence, without loss of generality, $B$ is an 
infinite set and
$t.d.(A)<\infty$. Let $T=K\otimes _kA={\begin{array}[t]{cl}
\lim\\
^{\stackrel{\rightarrow}{\scriptstyle E finite,\ E\subseteq 
B}}
\end{array}}T_{E}$, where $
T_{E}:=k_{E}\otimes _kA\subseteq T$ and $k_E:=k(E)$. Let us 
point out that,
for any
finite subset $E$ of $B$
and any prime ideal $P_{E}$ of $T_{E}$, $P_{E}T$ is a prime 
ideal of $T$.
Indeed, let $E$ be a finite subset of $B$ and $P_{E}$ a 
prime ideal of
$T_{E}$.
Then $T\cong K\otimes _{k_{E}}T_{E}$ and $P_{E}T=P_{E}
(K\otimes
_{k_{E}}T_{E})=K\otimes _{k_{E}}P_{E}$. Thus $\frac T{P_{E}
T}\cong
\frac{K\otimes
_{k_{E}}T_{E}}{K\otimes _{k_{E}}P_{E}}\cong K\otimes
_{k_{E}}\frac{T_{E}}{P_{E}}$. Note
that if $F$ is a field, $L=F(X_1,...,X_n)$ and $D$ is a 
domain containing
$F$,
then $L\otimes_FD$ ($\cong D[X_1,...,X_n]_{F[X_1,...,X_n]
\setminus\{0\}}$)
is a domain.
It follows that $\frac T{P_{E}T}$ is an integral domain, as 
desired, since
it is
a directed union of the domains $
k_{F}\otimes _{k_{E}}\frac{T_{E}}{P_{E}}$, where $F$ is a 
finite subset of
$B$
containing $E$.

Let $P\in \Spec(T)$ and $P_{E}:=P\cap T_{E},$ for each
finite subset $E$ of $B$. We claim that there exists a 
finite subset $E$ of
$B$
such that $P=P_{E}T$. Suppose by way of contradiction that 
for each finite
subset $E$ of $B$ we
have $P_{E}T\subset P$. Let $F$ be a finite subset of $B$. 
Assume that $%
P_{E}T=P_{F}T$ for each finite subset $E$ of $B$ that 
contains $F$. Let
$x\in P$.
Since $x\in T={\begin{array}[t]{ll}
\lim\\
^{^{\scriptstyle \rightarrow}}
\end{array}}T_{E},$ there exists a
finite subset $E_1$ of $B$ such that $x\in T_{E_1}$. Then 
$x\in P_{E_1}T.$
Thus $x\in P_{E_1\cup F}T=P_{F}T$. It follows that $P=P_{F}
T$, a
contradiction. Consequently, there exists a finite subset 
$E$ of $B$ such
that $F\subset E$ and $P_{F}T\subset P_{E}T$. Hence, by 
iterating the above
argument, we can construct an infinite
chain of prime ideals $P_{E_1}T\subset P_{E_2}
T\subset ...\subset
P_{E_n}T\subset ...\subset P,$ where the $E_j$ are finite 
subsets of $B$.
This is a
contradiction, since, by Proposition 4.1, $T$ is LFD. 
Therefore there exists
a finite subset $E$
of $B$ such that $P=P_{E}T$, proving the claim.

Let $P\subset Q$ be a chain of prime ideals of $T$. Then 
there exists a
common
finite subset $E$ of $B$ such that $P=P_{E}T$
and $Q=Q_{E}T$. We claim (*): $P\subset Q$ is saturated in 
$\Spec(T)$ if and
only
if $P_{E}\subset Q_{E}$ is saturated for each finite subset 
$E$ of $B$
such that $P=P_{E}T$ and $%
Q=Q_{E}T$. Indeed, assume that $P\subset Q$ is saturated 
and consider a
finite subset $E$ of $B$ such that $P=P_{E}T$ and $Q=Q_{E}
T$. Let $J$ be a
prime ideal of $T_{E}$ such that $P_{E}\subseteq J\subseteq 
Q_{E}$. Then $%
P_{E}T=P\subseteq JT\subseteq Q_{E}T=Q$. Since $\htt(\frac QP)
=1$ and $JT$ is
a
prime ideal of $T$, we obtain that either $JT=P=P_{E}T$ or 
$JT=Q=Q_{E}T$.
Since $T_{E}\subset _{\rightarrow }T$ is a faithfully flat 
homomorphism,
we conclude that either $J=P_{E}$ or $J=Q_{E}$
(see condition (i) in [2,
Exercise 16, p. 45]). Then $P_{E}\subset Q_{E}$ is 
saturated. Conversely,
suppose that $P_{E}\subset Q_{E}$ is saturated for each 
finite subset
$E$ of $B$
such that $P=P_{E}T$ and $Q=Q_{E}T$. Let $P^{\prime }$ be a 
prime ideal
of $T$
such that $P\subseteq P^{\prime }\subseteq Q$. There exists 
a finite
subset $F$ of $B$ satisfying $P=P_{F}T,P^{\prime }=P^
{\prime }_{F}T$
and $Q=Q_{F}T$.
Then $P_{F}\subseteq P^{\prime }_{F}\subseteq Q_{F}$. By 
hypothesis, $%
P_{F}\subset Q_{F}$ is saturated, so either $P^{\prime }_{F}
=P_{F}$ or $%
P^{\prime }_{F}=Q_{F}$. Hence, either $P^{\prime }=P$ or $P^
{\prime }=Q$.
Then $P\subset Q$ is saturated. This establishes the claim.

Now, assume that $A$ is a stably strong S-ring and let 
$P\subset Q$ be a
saturated chain in $\Spec(T)$. Then $P_{E}\subset Q_{E}$ is 
saturated for
each
finite subset $E$ of $B$ such that $P=P_{E}T$ and $Q=Q_{E}
T$. Hence
$P_{E}[X]\subset Q_{E}[X]$ is saturated, for each finite 
subset $E$ of
$B$ such that $P=P_{E}T$ and $Q=Q_{E}T$. We have $T[X]=
(K\otimes _kA)[X]=K\otimes _k(A[X])\cong {\begin{array}[t]
{cl}
\lim\\
^{\stackrel {\rightarrow}{E finite,\ E\subseteq B}}
\end{array}}(T_{E}[X])$. In view of
the equivalence (*), replacing $T$ by $T[X]$, $P$ by $P[X]$ 
and $Q$ by
$Q[X]$,
we conclude that $P[X]\subset Q[X]$ is saturated.
Therefore $T$ is a strong S-ring. Let $n\geq 1$ be an 
integer. Since $%
T[X_1,...,X_n]\cong K\otimes _k(A[X_1,...,X_n])$ and $A
[X_1,...,X_n]$ is
a stably strong S-ring, by repeating the earlier argument 
with $A$ replaced
by
$A[X_1,...,X_n]$, we can show that $K\otimes _k(A
[X_1,...,X_n])\cong
T[X_1,...,X_n]$
is a strong S-ring. Hence $T$ is a stably strong S-ring.

Now, suppose that $A$ is universally catenarian. We first 
recall (use
$E:=\emptyset$ in an earlier part of the proof)
that $\frac {K\otimes_kA}{K\otimes_kp}\cong K\otimes_k\frac 
Ap$ is a
domain, for any prime ideal $p$ of $A$. Furthermore, as $A 
\subset T$
satisfies GD,
one can easily check that $\Min(T) = \{K\otimes_kp : p \in 
\Min(A)\}$. It
follows that
$K\otimes_kA$ satisfies MPC, since $A$ satisfies MPC by 
hypothesis.
Moreover,
$T$ is LFD by Proposition 4.1. Let $P\subset Q$
be a saturated chain of prime ideals of $T$. Then $P_{E}
\subset
Q_{E}$ is saturated for each finite subset $E$ of $B$ such 
that $P=P_{E}T$
and $Q=Q_{E}T$. Take a finite subset $E=\{X_1,...,X_n\}$ of 
$B$ and set $S_E=k[X_1,...,X_n]\setminus\{0\}$. Then
$T_{E}\cong S^{-1}_{E}A[X_1,...,X_n]$ is (universally)
catenarian,
by the hypothesis on $A$. Hence, $\htt(Q_{E})=1+\htt(P_{E})$ 
for each finite
subset
$E$ of $B$ such that $P=P_{E}T$ and $Q=Q_{E}T$. On the 
other hand, we claim
that $\htt(P)=\sup\{\htt(P_{E}):E$ is a finite subset of $B$ 
such that
$P=P_{E}T\}$
and $%
\htt(Q)=\sup\{\htt(Q_{E}):E$ is a finite subset of $B$ such that 
$Q=Q_{E}T\}$.

Indeed, let $E$ be a finite subset of $B$ such that $P=P_{E}
T$. Since the
homomorphism $T_{E}\subset_{\rightarrow }T$ satisfies GD, 
we have $%
\htt(P_{E})\leq \htt(P)$. Hence $\sup\{\htt(P_{E}):E$ is a finite 
subset of $B$
such
that $P=P_{E}T\}\leq \htt(P)$. Since $T$ is LFD, $\htt(P)$ is
finite. Let $
P_0\subset P_1\subset ...\subset P_h=P$ be a chain of prime 
ideals of $T$
such that
$h=\htt(P)$. There exists a common finite subset $E$ of $B$ 
such that $
P_i=P_{iE}T$, for $i=0,...,h$. Then $P_{0E}\subset P_{1E}
\subset ...\subset
P_{hE}$ is a chain of distinct prime ideals in $T_{E}$, 
since the
homomorphism
$T_{E}\rightarrow T$ is faithfully flat. Hence $h=\htt(P)\leq
\htt(P_{hE})=\htt(P_{E})$. It follows that $\htt(P)\leq \sup\{\htt(P_
{E}):E$ is
a finite subset of $B$ such that $P=P_{E}T\}$. This 
establishes the
above claim. We conclude that $\htt(Q)=1+\htt(P)$. Hence $T$ is 
catenarian.
Since $
T[X_1,...,X_n]\cong K\otimes _k(A[X_1,...,X_n])$, an
argument similar to the above,
with $A$ replaced by $A[X_1,...,X_n]$, shows that $T$ is
universally catenarian and the proof is complete. 
$\diamondsuit$\bigskip

\noindent {\bf Proof of Theorem 4.13.} We have
\(K\otimes _kA \cong K\otimes_{k(B)}(k(B)\otimes _kA)
 \cong \linebreak K\otimes _{L}(L\otimes
_{k(B)}(k(B)\otimes_kA))\).
Since $[L:k(B)]<\infty$, we have $K=k(B)(x_1,...,x_n)$ for 
some
$x_1,...,x_n\in L$. So $L\cong \frac {k(B)[X_1,...,X_n]}I$, 
for some
prime ideal $I$ of $k(B)[X_1,...,X_n]$. It follows that
$L\otimes_{k(B)}(k(B)\otimes_kA)\cong \frac
{(k(B)\otimes_kA)[X_1,...,X_n]}J$,
where $J=I\otimes_{k(B)}(k(B)\otimes_kA)$. By Proposition 4.14, 
$k(B)\otimes _kA$
is a stably strong
S-ring (resp., universally catenarian) if $A$ is. Thus, if 
$A$ is a
stably
strong S-ring (
resp., universally catenarian), $L\otimes_{k(B)}(k(B)
\otimes_kA)$ is so (we
have just used the easy fact that the class of stably 
strong S-rings is
closed under formation of factor rings).
Then,
by Proposition 4.5, the result follows, since $K$ is a 
purely inseparable
extension of $L$. $\diamondsuit$
\end{section}\bigskip

\begin{section}{Examples}

This section displays some examples showing that several 
results of
Section 4 concerning the strong S-property and catenarity
of $K\otimes_kA$ fail, in general, when the field extension 
$K$
is no longer algebraic over $k$. Our last example, Example 
5.5, shows
clearly that the study of the spectrum of $A\otimes_kB$ 
becomes more
intricate
if one moves beyond the context where at least one of $A$, 
$B$ is a field
extension of $k$.\bigskip

In order to provide some background for the present section,
we recall the following definitions and results from [26]. 
A domain $A$ is
called an AF-domain if $A$ is a
$k$-algebra of finite transcendence degree over $k$
such that $\htt(p)+ t.d.(\frac Ap:k) = t.d.(A:k)$ for each 
$p\in \Spec(A)$.
Finitely generated $k$-algebras (that are domains) and field
extensions
of finite transcendence degree over $k$ are AF-domains.
Let $A$ be a k-algebra, $p$
a prime ideal of $A$ and $0\leq d\leq s$ be integers. Set
\[ \triangle(s,d,p):=\htt(p[X_1,...,X_s])+ \min(s,d+t.d.(\frac 
Ap:k)),\]
\[D(s,d,A):=\max\{\triangle (s,d,p):p\in \Spec(A)\}.\]
Wadsworth's main two results
relative to the Krull dimension of tensor products of AF-
domains read as
follows.
If $A$ is an AF-domain and $R$ is any $k$-algebra, then
$\dim (A\otimes _kR)=D(t.d.(A:k),\dim (A),R)$ [26, Theorem 
3.7].
If, in addition, $R$ is an AF-domain, then
$\dim (A\otimes_kR)=\min(\dim (A)+t.d.(R:k),t.d.(A:k)+\dim
(R))$ [26, Theorem
3.8].\bigskip

We turn now to our examples. It is still an open problem to 
know whether
$K\otimes_kA$ is a strong
S-ring (resp., catenarian) when $K$ is an algebraic field 
extension of
$k$ and $A$ is a strong S-ring (resp., catenarian such that 
$K\otimes_kA$
satisfies MPC). However, for the case where $K$ is a 
transcendental
field extension of $k$, the answer is negative, as
illustrated by the following two examples.
\begin{example} Let $k$ be a field. There
exists a strong S-domain $A$ that is a $k$-algebra such 
that $L\otimes _kA$
is a strong
S-ring for any algebraic field extension $L$ of $k$, while 
$K\otimes _kA$
is not a strong S-ring for some transcendental field 
extension $K$ of
$k$.\end{example}

Our example draws on [8, Example 3],
which we assume that the reader has at hand. Let $k$ be a 
field and
$k^{\prime }$ an algebraic closure of $k$. Let $(V_1, M^
{\prime}_1)$ be the
valuation
domain of the $Y_3$-adic valuation on $k^{\prime}(Y_1,Y_2)
[Y_3]$ . Let $V^*$
be a discrete rank-one
valuation domain of $k^{\prime}(Y_1,Y_2)$ of the form $k^
{\prime}+N$
and let $V$ be the pullback $\varphi^{-1}(V^*)$, where
$\varphi:V_1\rightarrow k^{\prime}(Y_1,Y_2)$ is the 
canonical
homomorphism. It is easily seen that $V$ is a rank-two 
valuation domain of
the form $k^{\prime}+M_1$. Moreover, if $p_1$ is the height 
1 prime
ideal of $V$, $V_{p_1}=V_1$. Finally, let $W$ be the 
valuation domain of the
$(Y_3+1)$-adic valuation
on $k^{\prime}(Y_1,Y_2)[Y_3]$. Then $W$ is a DVR of the form
$k'(Y_1,Y_2)+M_2$. Set $A=k^{\prime }+M$, where
$M=M_1\cap M_2$. It is shown in [8, Example 3] that $A$ is 
a two-dimensional
local strong S-domain
with the following features: $\dim A[X,Y]=5$ (hence $A[X]$ 
is not a strong
S-domain),
the quotient field of $A$ is $k^{\prime}(Y_1,Y_2,Y_3)$, and
the prime ideals of $A$ are $(0)\subset p\subset M$ with 
$A_p=V_1$.
By Proposition 4.6, $L\otimes _kA$ is a strong S-ring, for 
any algebraic
field extension $L$ of $k$.
On the other hand, by [26, Theorem 3.7], $\dim ((k(X)
\otimes _kA)[Y])=\dim
(k(X)[Y]\otimes _kA)=D(2,1,A)$, since
$k(X)[Y]$ is an AF-domain. We have
\begin{eqnarray*}
\triangle (2,1,(0)) &=& \min(2,1+t.d.(A:k))\\
 &=& \min(2,4)=2. \\
\triangle (2,1,p) &=& \htt(p[X,Y])+ \min(2,1+t.d.(\frac 
Ap:k))\\
 &=& \htt(pA_p[X,Y])+\min(2,1+t.d.(\frac
{A_p}{pA_p}:k))\\
&=& \htt(pA_p)+\min(2,1+t.d.(\frac{V_1}{M^{\prime}_1}:k))\ 
\hbox{(since
$A_p$ is a DVR)}\\
&=& 1+\min(2,3)=3. \\
\triangle (2,1,M) &=& \htt(M[X,Y])+\min(2,1+t.d.(\frac 
AM:k)) \\
 &=& \dim A[X,Y]-2+\min(2,1)= 4.
\end{eqnarray*}
Hence $\dim (k(X)\otimes _kA)[Y])=4$. Furthermore, $\dim
(k(X)\otimes_kA)=D(1,0,A)$. We have
\begin{eqnarray*}
\triangle (1,0,(0)) &=& \min(1,t.d.(A:k))\\
 &=& \min(1,3)=1. \\
\triangle (1,0,p) &=& \htt(p[X])+\min(1,t.d.(\frac Ap:k)) 
\\
 &=& \htt(pA_p[X])+\min(1,2)\\
 &=& \htt(pA_p)+ 1=2\ \hbox{(since $A_p$ is a DVR)}.\\
\triangle (1,0,M) &=& \htt(M[X])+\min(1,t.d.(\frac AM:k)) 
\\
 &=& \htt(M)+\min(1,0)=2\ \hbox{(since 
$A$ is a strong
S-domain).}
\end{eqnarray*}
Hence, $\dim(k(X))\otimes_kA)=2$.
Consequently, $\dim((k(X)\otimes _kA)[Y])=4\neq1+2= 1+\dim(k
(X)\otimes
_kA)$. Let $K=k(X)$.
Therefore, by [18, Theorem 39], $K\otimes _kA$ is not a 
strong S-ring.
$\diamondsuit$
\begin{example} Let $k$ be a field. There
exists a catenarian domain $A$ that is a $k$-algebra such 
that $L\otimes
_kA$ is
catenarian for any algebraic field extension $L$ of $k$, 
while $K\otimes
_kA$
is not catenarian for some transcendental field extension 
$K$ of $k$.\end{example}

Let $k$ be a field and $k^{\prime }$ an algebraic closure 
of $k$.
Let $V:=k^{\prime }(X_1,X_2)[Y]_{(Y)}=k^{\prime}(X_1,X_2)
+m$, where $m:=YV$.
Let $A:=k^{\prime }(X_1)+m$. Clearly, $A$
is catenarian while $A[Z]$ is
not catenarian, as the following chains of prime
ideals of $A[Z]$ are saturated:

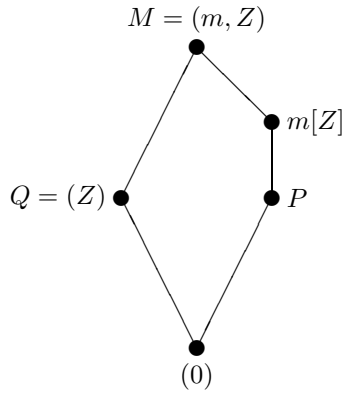
\[\setlength{\unitlength}{1mm}
\begin{picture}(40,50)(10,-5)
\put(20,0){\line(-1,2){10}}
\put(10,20){\line(1,2){10}}
\put(20,40){\line(1,-1){10}}
\put(30,30){\line(0,-1){10}}
\put(30,20){\line(-1,-2){10}}
\put(20,40){\circle*{2}}
\put(20,42){\makebox(0,0)[b]{$M=(m,Z)$}}
\put(10,20){\circle*{2}}
\put(8,20){\makebox(0,0)[r]{$Q=(Z)$}}
\put(20,0){\circle*{2}}
\put(20,-2){\makebox(0,0)[t]{$(0)$}}
\put(30,20){\circle*{2}}
\put(32,20){\makebox(0,0)[l]{$P$}}
\put(30,30){\circle*{2}}
\put(32,30){\makebox(0,0)[l]{$m[Z]$}}
\end{picture}\]
\noindent where $P$ is an upper to $(0)$ (cf. [8, Example 
5]).
By Proposition 4.6, $L\otimes _kA$ is catenarian for any
algebraic field extension $L$ of $k.$ On the other hand,
$\frac{S^{-1}A[Z]}{%
S^{-1}m[Z]}$ $\cong \frac{k(Z)\otimes _kA}{k(Z)\otimes _km}
\cong k(Z)\otimes
_k\frac Am\cong k(Z)\otimes _kk^{\prime }(X_1)$; let $S=k
[Z]\setminus
\{0\}$. Therefore
$\dim (\frac{S^{-1}A[Z]}{S^{-1}m[Z]})=1$ by [23, Theorem 
3.1]. Hence
$S^{-1}m[Z]$ is not a
maximal ideal of $S^{-1}A[Z]$, whence
there exists an upper $M_1$ to $m$ such that $M_1\cap 
S=\emptyset$. By
[9, Theorem B, p. 167],
$l(M)=l(M_1)$, where $l(M)$ (resp., $l(M_1)$) denotes the 
set of lengths
of saturated chains of prime ideals between $(0)$ and $M$ 
(resp., $M_1$).
Then there exist two
saturated chains of prime ideals in $A[Z]$ of the form:
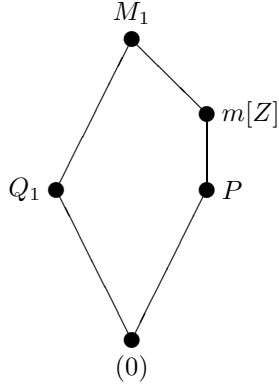
\[\setlength{\unitlength}{1mm}
\begin{picture}(40,50)(10,-5)
\put(20,0){\line(-1,2){10}}
\put(10,20){\line(1,2){10}}
\put(20,40){\line(1,-1){10}}
\put(30,30){\line(0,-1){10}}
\put(30,20){\line(-1,-2){10}}
\put(20,40){\circle*{2}}
\put(20,42){\makebox(0,0)[b]{$M_1$}}
\put(10,20){\circle*{2}}
\put(8,20){\makebox(0,0)[r]{$Q_1$}}
\put(20,0){\circle*{2}}
\put(20,-2){\makebox(0,0)[t]{$(0)$}}
\put(30,20){\circle*{2}}
\put(32,20){\makebox(0,0)[l]{$P$}}
\put(30,30){\circle*{2}}
\put(32,30){\makebox(0,0)[l]{$m[Z]$}}
\end{picture}\]

\noindent where $Q_1$ in an upper to $(0)$.
Consequently, $K\otimes _kA\cong S^{-1}A[Z]$ is not 
catenarian, where
$K:=k(Z)$. $\diamondsuit$\bigskip

The next two examples show that Proposition 4.4 fails in 
general when
$K$ is no longer algebraic over $k$.
\begin{example} There exists a $k$-algebra $A$ 
which is not an
S-domain and a field extension $K$ of $k$ such that $1\leq 
t.d.(K:k)<
\infty $ and $K\otimes _kA$ is a strong S-ring.\end{example}

Let $V:=k(X)[Y]_{(Y)}=k(X)+m$, where $m:=YV$, and let 
$A:=k+m$. We have $
\htt(m)=1$ and $\htt(m[Z])=\htt(m[Z,T])=2$ [8, Example 5]. Thus, 
$A$ is not an
S-domain.
Let $K:=k(Z)$. We claim that $K\otimes _kA\cong S^{-1}A[Z]$ 
is a strong
S-domain, where $S:=k[Z]\setminus\{0\}$. Notice first that 
$S^{-1}m[Z]$ is a
maximal
ideal of $S^{-1}A[Z]$, as $\frac{S^{-1}A[Z]}{S^{-1}m[Z]}
\cong
\frac{k(Z)\otimes _kA}{k(Z)\otimes _km}\cong k(Z)\otimes 
_k\frac Am\cong
k(Z)$.
Now, let $P\subset Q$ be a
pair of adjacent prime ideals of $A[Z]$ that are disjoint 
from $S$. Two
cases
are possible. If $P=(0)$, then $\htt(Q)=1$. Since
$k(Z)\otimes _kA\cong S^{-1}A[Z]$ is an S-domain, $\htt(Q[T])
=1$. If
$P$ is an upper to $(0)$, $Q$ necessarily contracts to $m$ 
in A
and hence $Q=m[Z]$, since $Q\cap S=\emptyset$ and $S^{-1}m
[Z]\in
\Max(S^{-1}A[Z])$. Therefore
$(0)\subset P\subset m[Z]=Q$ is a
saturated chain in $\Spec(A[Z])$. Then $(0)\subset
P[T]\subset m[Z,T]=Q[T]$ is a saturated chain in $\Spec(A
[Z,T])$.
Consequently, in
both
cases, $P[T]\subset Q[T]$ is saturated. It follows that 
$K\otimes _kA\cong
S^{-1}A[Z]$
is a strong S-domain, as desired. $\diamondsuit$
\begin{example} There exists a $k$-algebra 
$A$ which is not
a catenarian domain and a field extension $K$ of $k$ such 
that $1\leq
t.d.(K:k)<
\infty $ and $K\otimes _kA$ is catenarian.\end{example}

Let $V:=k(X)[Y]_{(Y)}=k(X)+m,$ where $m:=YV.$ Let $R:=k+m$. 
Clearly, $R$ is
a
one-dimensional integrally closed domain. There exist two 
saturated
chains of prime ideals of $R[Z]$, as in Example 5.2, of the 
form:
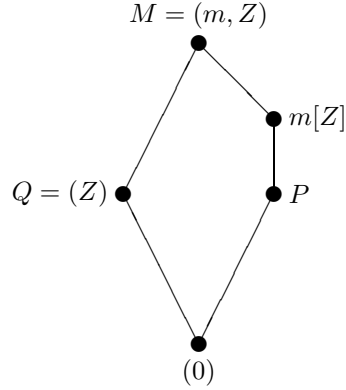
\[\setlength{\unitlength}{1mm}
\begin{picture}(40,50)(10,-5)
\put(20,0){\line(-1,2){10}}
\put(10,20){\line(1,2){10}}
\put(20,40){\line(1,-1){10}}
\put(30,30){\line(0,-1){10}}
\put(30,20){\line(-1,-2){10}}
\put(20,40){\circle*{2}}
\put(20,42){\makebox(0,0)[b]{$M=(m,Z)$}}
\put(10,20){\circle*{2}}
\put(8,20){\makebox(0,0)[r]{$Q=(Z)$}}
\put(20,0){\circle*{2}}
\put(20,-2){\makebox(0,0)[t]{$(0)$}}
\put(30,20){\circle*{2}}
\put(32,20){\makebox(0,0)[l]{$P$}}
\put(30,30){\circle*{2}}
\put(32,30){\makebox(0,0)[l]{$m[Z]$}}
\end{picture}\]
\noindent Let $A:=R[Z]$. Then $A$ is not catenarian. We 
next prove
that $K\otimes _kA\cong S^{-1}R[Z,T]$ is catenarian, where 
$K:=k(T)$
and $S:=k[T]\setminus\{0\}$.

Notice first
that $\htt(m[Z,T])=2$ [8, Example 5]. Further, one may easily
check, via [26, Theorem 3.7], that $\dim (K\otimes _kA)=
\dim (k(T)[Z]\otimes_kR)=D(2,1,R)= 3$, since $k(T)[Z]$ is 
an AF-domain.
Now, let $P \subseteq Q$ be a pair of prime ideals of $R
[Z,T]$
such that $Q \cap S= \emptyset$. We claim that
$\htt(Q) = \htt(P) + \htt(Q/P)$. Without loss of generality, we 
may assume that
$\htt(Q)=3$. Necessarily, $Q$ contracts to $m$ in $R=k+m$. 
Moreover,
$Q$ cannot be an upper to an upper to $m$ in $R[Z,T]$; 
otherwise $\htt(Q)=4$.
Hence, either $Q=M_1[T]$ or $Q=M_2[Z]$, where $M_1$ is an 
upper to $m$
in $R[Z]$ and $M_2$ is an upper to $m$ in $R[T]$. Assume 
that $Q=M[T]$,
where $M$ is
an upper to $m$ in $R[Z]$. In case $P\cap R=m$, we are 
done, since here
$P=m[Z,T]$.
We may then assume that $P\cap R=(0)$. Three cases are 
possible. If
$P$ is an upper to an upper to $(0)$ in $R[Z,T]$, then $ht
(P)=2$, and we are
done.
If $P=P_1[Z]$, where $P_1$ is an upper to $(0)$ in $R[T]$, 
then $P\cap
R[T]=P_1\subset Q\cap R[T]=(M\cap R)[T]=m[T]$. Hence $P=P_1
[Z]\subset
m[Z,T]$. Thus $\htt(\frac {Q}{P})=2$ and $\htt(P)=1$, as 
desired. Assume now that
$P=P_2[T]$, where $P_2$ is an upper to $(0)$ in $R[Z]$. We
have $Q=M[T]$ is an upper to $m[T]$ in $(R[T])[Z]$ and $P$ 
is
an upper to $(0)$ in $(R[T])[Z]$. If $\htt(\frac {Q}{P})
=1<\htt((m[T])[Z])=2$, then by [9, Proposition 2.2] and
[14, Proposition 1.1, p. 742], $P\subset m[Z,T]$ (since
$R[T]$ is integrally closed), a contradiction. Thus, $ht
(\frac {Q}{P})=2$
and $\htt(P)=1$, as desired.
A similar argument applies to the case where $Q=M[Z]$,
where $M$ is an upper to $m$ in $R[T]$. Consequently, 
$K\otimes_kA$ is
catenarian. $\diamondsuit$\bigskip

To emphasize the importance of $K$ being a field in Theorem 
4.13, we close
this section with an example of two discrete rank-one 
valuation
domains,
hence universally catenarian, the tensor product of which 
is not
catenarian.\bigskip

\noindent {\bf Example 5.5.} There exists a discrete rank-one valuation
domain $V$
such that $t.d.(V:k)< \infty $ and $V\otimes _kV$ is not 
catenarian.\bigskip

Consider the $k$-algebra homomorphism $\varphi
:k[X,Y]\rightarrow k[[t]]$ such that\linebreak $\varphi (X)=t$ and 
$\varphi (Y)=s:=
{\Sigma}_{n\geq 1}t^{n!}$. Since $s$ is known to be 
transcendental over
$k(t)$,
$\varphi$ is injective. This induces an embedding $\overline
{\varphi }
:k(X,Y)\rightarrow k((t))$ of fields. Put $V=\overline
{\varphi}
^{-1}(k[[t]])$. It is easy to check that $V$ is a discrete 
rank-one
valuation overring of $k[X,Y]$ of the form $k + m$, where 
$m:=XV$. For
convenience, put $A=B:=V$. We have $\dim(A\otimes_kB)=\dim 
(V\otimes
_kV)=\dim (V)+t.d.(V:k)=1+2=3$ [26,
Corollary 4.2] and $\htt(m\otimes_kV)=\htt(m[X,Y])=\htt(m)=1$ [4,
Lemma
1.4]. Since $\htt(\frac{m\otimes _kV+V\otimes _km}{m\otimes 
_kV})\leq \dim (%
\frac{V\otimes _kV}{m\otimes _kV})=\dim (V)=1$, we obtain 
$\htt(\frac{m\otimes
_kV+V\otimes _km}{m\otimes _kV})=1$. On the other hand, in 
view of
[26, Proposition 2.3], the height of no prime ideal of 
$A\otimes_kB$
contracting to $(0)$ in $A$ and to $(0)$ in $B$ can reach
$\dim(A\otimes_kB)=3$,
since $\dim (k(X,Y)\otimes
_kk(X,Y))=2$. Therefore, $\htt(m\otimes _kV+V\otimes _km)=3$. 
Hence
$\Spec(V\otimes_kV)$ contains the following two saturated 
chains:
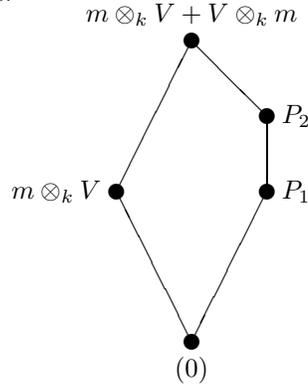
\[\setlength{\unitlength}{1mm}
\begin{picture}(40,50)(10,-5)
\put(20,0){\line(-1,2){10}}
\put(10,20){\line(1,2){10}}
\put(20,40){\line(1,-1){10}}
\put(30,30){\line(0,-1){10}}
\put(30,20){\line(-1,-2){10}}
\put(20,40){\circle*{2}}
\put(20,42){\makebox(0,0)[b]{$m\otimes_k V + V\otimes_k m$}}
\put(10,20){\circle*{2}}
\put(8,20){\makebox(0,0)[r]{$m\otimes_k V$}}
\put(20,0){\circle*{2}}
\put(20,-2){\makebox(0,0)[t]{$(0)$}}
\put(30,20){\circle*{2}}
\put(32,20){\makebox(0,0)[l]{$P_1$}}
\put(30,30){\circle*{2}}
\put(32,30){\makebox(0,0)[l]{$P_2$}}
\end{picture}\]
\noindent where $P_i\cap A=P_i\cap B=(0)$, for $i=1,2$. 
Consequently,
$V\otimes _kV$ is
not catenarian. $\diamondsuit$
\end{section}\bigskip\bigskip

\noindent{\bf References}

\begin{list}{}{\topsep=1mm \itemsep=0mm \parsep=0mm 
\labelsep=2mm
\labelwidth=14mm}

\item [{[1]}] D.F. Anderson, A. Bouvier, D.E. Dobbs, M.
Fontana, S. Kabbaj, On Jaffard domains, Expo. Math.  6  (1988) 145-175.\medskip

\item [{[2]}] M.F Atiyah, I.G. Macdonald, Introduction to
commutative algebra, Addi-son-Wesley, Reading, Mass., 1969.\medskip

\item [{[3]}] S. Bouchiba, F. Girolami, S. Kabbaj, The
dimension of tensor products of AF-rings, Lecture Notes 
Pure Appl.
Math., Dekker, New York,  185  (1997) 141-154.\medskip

\item [{[4]}] S. Bouchiba, F. Girolami, S. Kabbaj, The
dimension of tensor products of algebras arising from 
pullbacks, J.
Pure Appl. Algebra  137  (1999) 125-138.\medskip

\item [{[5]}] A. Bouvier, D.E. Dobbs, M. Fontana, Universally
catenarian integral domains, Advances in Math.  72  (1988) 211-238.\medskip

\item [{[6]}] A. Bouvier, D.E. Dobbs, M. Fontana, Two sufficient
conditions for universal catenarity, Comm. Algebra  15  (1987) 861-872.\medskip

\item [{[7]}] A. Bouvier, M. Fontana, The catenarian 
property
of the polynomial rings over a Pr\"ufer domain, in S\'em. 
Alg\`ebre P.
Dubreil
et M.P. Malliavin, Lecture Notes in Math.,  1146,
Springer-Verlag, Berlin-New York, 1985, 340-354.\medskip
  
\item [{[8]}] J.W. Brewer, P.R. Montgomery, E.A. Rutter, 
W.J. Heinzer, Krull dimension of polynomial rings, Lecture Notes in 
Math.,  311, Springer-Verlag, Berlin-New York, 1972, 26-45.\medskip

\item [{[9]}] A.M. De Souza Doering, Y. Lequain,
Chains of
prime ideals in polynomial rings, J. Algebra  78  (1982) 163-180.\medskip

\item [{[10]}] D.E. Dobbs, M. Fontana, S. Kabbaj, Direct
limits of Jaffard domains and S-domains, Comm. Math. Univ. 
St. Pauli  39  (1990) 143-155.\medskip

\item [{[11]}] M. Fontana, S. Kabbaj, On the Krull 
and valuative
dimension of $D+XD_S[X]$ domains, J. Pure Appl. 
Algebra  63  (1990) 231-245.\medskip

\item [{[12]}] R. Gilmer, Multiplicative ideal 
theory, Dekker, New
York, 1972.\medskip

\item [{[13]}] A. Grothendieck, J.A. Dieudonn\'e, 
El\'ements
de g\'eom\'etrie alg\'ebrique, Springer-Verlag, Berlin, 1971.\medskip

\item [{[14]}] E.G. Houston, S. Mcadam, Chains of 
primes in
Noetherian rings, Indiana Univ. Math. J.  24  (1975) 741-753.\medskip

\item [{[15]}] P. Jaffard, Th\'eorie de la dimension 
dans les anneaux
de polyn\^omes, M\'em. Sc. Math.,  146 , Gauthier-Villars, 
Paris, 1960.\medskip

\item [{[16]}] S. Kabbaj, La formule de la dimension 
pour les
S-domaines
forts universels, Boll. Un. Mat. Ital.  D-VI 5  (1) (1986) 145-161.\medskip

\item [{[17]}] S. Kabbaj, Sur les S-domaines forts 
de Kaplansky,
J. Algebra  137  (2) (1991) 400-415.\medskip

\item [{[18]}] I. Kaplansky, Commutative rings, 
rev. ed., University of 
Chicago Press, Chicago, 1974.\medskip

\item [{[19]}] S. Malik, J.L. Mott, Strong S-
domains, J. Pure
Appl. Algebra  28  (1983) 249-264.\medskip

\item [{[20]}] H. Matsumura, Commutative ring 
theory, Cambridge
University Press, Cambridge, 1989.\medskip

\item [{[21]}] M. Nagata, Local rings, 
Interscience, New York, 1962.\medskip

\item [{[22]}] L.J. Ratliff, On quasi-unmixed local 
domains, the
altitude formula, and the chain condition for prime ideals, 
II,
Amer. J. Math.  92  (1970) 99-144.\medskip

\item [{[23]}] R.Y. Sharp, The dimension of the 
tensor product
of two field extensions, Bull. London Math. Soc.  9  (1977) 42-48.\medskip

\item [{[24]}] R.Y. Sharp, P. Vamos, The dimension 
of the tensor
product of a finite number of field extensions, J. Pure 
Appl. Algebra  10  (1977) 249-252.\medskip

\item [{[25]}] P. Vamos, On the minimal prime ideals 
of a tensor
product of two fields, Math. Proc. Camb. Phil. Soc.  84  (1978) 25-35.\medskip

\item [{[26]}] A.R. Wadsworth, The Krull dimension 
of tensor products
of commutative algebras over a field, J. London Math. 
Soc.  19  (1979) 391-401.\medskip

\item [{[27]}] O. Zariski, P. Samuel, Commutative 
algebra,
Vol. I, Van Nostrand, Princeton, 1958.\medskip

\item [{[28]}] O. Zariski, P. Samuel, Commutative 
algebra,
Vol. II, Van Nostrand, Princeton, 1960.
\end{list}
\end{document}